\documentclass[11pt,leqno]{amsart}
\usepackage{times,amsmath,amsthm,amssymb,latexsym,amscd}
\usepackage[all]{xy}
\usepackage{graphicx}

\newtheoremstyle{fancy}{}{}{\itshape}{}{\textsc\bgroup}{.\egroup}{ }{}
\newtheoremstyle{fanci}{}{}{\rm}{}{\textsc\bgroup}{.\egroup}{ }{}
\theoremstyle{fancy}
\newcounter{intro}

\newcounter{test}

\numberwithin{test}{subsection} 
\numberwithin{equation}{section} 
\newtheorem{cor}[test]{Corollary}
\newtheorem{lemma}[test]{Lemma}
\newtheorem{prop}[test]{Proposition}
\newtheorem{thm}[test]{Theorem}
\theoremstyle{fanci}
\newtheorem{dfn}[test]{Definition}
\newtheorem{ex}[test]{example}

\newcommand{\secref}[1]{Section~\ref{#1}}

\newcommand{\propref}[1]{Prop\-o\-si\-tion~\ref{#1}}
\newcommand{\corref}[1]{Cor\-ol\-lary~\ref{#1}}

\newcommand{\cref}[1]{Corollary~\ref{#1}}

\let\abs=\envert

\newcommand{\lqp}{L(q:p_1, \ldots, p_n)}

\newcommand{\iqn}{I_{0}(q,n)}
\newcommand{\itqn}{\widetilde{I}_{0}(q,n)}

\begin{document}

\title[Orbifold Lens Spaces That Are Isospectral But Not Isometric]{Orbifold Lens Spaces That Are Isospectral But Not Isometric}

\author[N. Shams]{Naveed Shams ul Bari}
\address{Naveed Shams ul Bari \\ 42/2 Lane 21, DHA Phase 7, Karachi, Pakistan - Director of Quantitative Methods at Center for Financial Training and Research (CFTR) in Karachi, Pakistan}
\email{bari.naveed@yahoo.com}


\begin{abstract}
We answer Mark Kac's famous question \cite{K}, ``can one
hear the shape of a drum?'' in the negative
for orbifolds that are spherical space forms. This is done by extending the
techniques developed by A. Ikeda on Lens Spaces to the orbifold setting. Several
results are proved to show that with certain restrictions on the
dimensionalities of orbifold Lens spaces we can obtain infinitely many pairs of
isospectral non-isometric Lens spaces. These results are then generalized to
show that for any dimension greater than 8 we can have pairs of isospectral
non-isometric orbifold Lens spaces. 

\end{abstract}

\maketitle

\tableofcontents

\section{Introduction}\label{introduction}

Given a closed Riemannian manifold $(M,g)$, the eigenvalue spectrum of the 
associated Laplace Beltrami operator will be referred to as the spectrum of 
$(M,g)$.  The inverse spectral problem asks the extent to which the spectrum 
encodes the geometry of $(M,g)$.  While various geometric invariants such as 
dimension, volume and total scalar curvature are spectrally determined, 
numerous examples of isospectral Riemannian manifolds, i.e., manifolds with 
the same spectrum, show that the spectrum does not fully encode the geometry.  
Not surprisingly, the earliest examples of isospectral manifolds were manifolds 
of constant curvature including flat tori (\cite{M}), hyperbolic manifolds (\cite{V}), 
and spherical space forms (\cite{I1}, \cite{I2} and \cite{Gi}).  In particular, Lens 
spaces are quotients of round spheres by cyclic groups of orthogonal transformations 
that act freely on the sphere.  Lens spaces have provided a rich source of isospectral 
manifolds with interesting properties.  In addition to the work of Ikeda cited above, 
see the recent results of Gornet and McGowan \cite{GoM}.

In this paper we generalize this theme to the category of Riemannian Orbifolds. 

A smooth \emph{orbifold} is a topological space that is locally modelled on an
orbit space of $\mathbf{R}^{n}$ under the action of a finite group of diffeomorphisms. 
\emph{Riemannian} orbifolds are spaces that are locally modelled on quotients of Riemannian 
manifolds by finite groups of isometries. Orbifolds have wide applicability, for example, in 
the study of 3-manifolds and in string theory.

The tools of spectral geometry can be transferred to the setting of Riemannian
orbifolds by using their well-behaved local structure (see \cite{Chi}, \cite{S1} \cite{S2}).
As in the manifold setting, the spectrum of the Laplace operator of a compact
Riemannian orbifold is a sequence $0 \leq \lambda_{1} \leq \lambda_{2} \leq \lambda_{3} \leq \ldots 
\uparrow \infty$ where each eigenvalue is repeated according to its finite multiplicity. 
We say that two orbifolds are isospectral if their Laplace spectra agree.

The literature on inverse spectral problems on orbifolds is less developed than that for manifolds.  
Examples of isospectral orbifolds include pairs with boundary (\cite{BCDS} and \cite{BW}); isospectral 
flat 2-orbifolds (\cite{DR}); arbitrarily large finite families of isospectral orbifolds (\cite{SSW}); 
isospectral orbifolds with different maximal isotropy orders (\cite{RSW}); and isospectral deformation 
of metrics on an orbifold quotient of a nilmanifold (\cite{PS}).

In this article, we study the spectrum of orbifold Lens spaces, i.e., quotients of round spheres by 
cyclic groups of orthogonal transformations that have fixed points on the sphere.  Generalizing the 
work of Ikeda (see \cite{I1}, \cite{I2} and \cite{IY}) we construct the generating function for the 
spectrum and systematically construct isospectral orbifold Lens spaces.  Section two introduces the 
orbifold Lens spaces and their generating functions.  In section 3, we will develop the proofs of our 
main theorems. We will first prove:

\smallskip

 \textbf{Theorem 3.1.6.} 
 
\begin{enumerate} 
\item [(i)] {\it Let $p \geq 5$ (alt. $p \geq 3$) be an odd prime and let $m
\geq 2$ (alt. $m \geq 3$) be any positive integer. Let $q = p^{m}$. Then there
exist at least two $(q-6)$-dimensional orbifold lens spaces with fundamental
groups of order $p^{m}$ 
which are isospectral but not isometric.}
 
\item[(ii)] {\it Let $p_{1},p_{2}$ be odd primes such that $q=p_{1} \cdot p_{2}
\geq 3$. Then there exist at least two $(q-6)$-dimensional orbifold 
lens spaces with fundamental groups of order $p_{1} \cdot p_{2}$ which are
isospectral but not isometric. }

\item[(iii)] {\it Let $q=2^{m}$ where $m \geq 6$ is any positive integer. Then
there exist at least two $(q-5)$-dimensional orbifold lens spaces 
with fundamental groups of order $2^{m}$ which are isospectral but not
isometric. }

\item[(iv)] {\it Let $q=2p$, where $p$ is an odd prime and $p \geq 7$. Then
there exist at least two $(q-5)$-dimensional orbifold lens spaces with
fundamental groups of order $2p$ which are isospectral but not isometric. }

\end{enumerate}

To prove these results we proceed as follows: 
\smallskip

\begin{enumerate}

\item Depending on the number of $p^{i}$ (alt. $p_{1},p_{2}$) divisors of
$q=p^{m}$ (alt. $q=p_{1} \cdot p_{2}$), we reformulate the 
generating function in terms of rational polynomial functions. 

\item Then we classify the number of generating functions that we will get by
imposing different conditions on the domain values of these 
polynomial functions. 

\item We prove sufficiency conditions on the number of generating functions that
would guarantee isospectrality for non-isometric orbifold 
lens spaces. 

\end{enumerate} 

The techniques used to prove these results parallel similar techniques from the
manifold lens space setting used in \cite{I1}. 

Generalizing this technique, we will get our second set of main results: 

\smallskip

\textbf{Theorem 3.2.5.} {\it Let $W \in \{0,1,2, \ldots \}$. }

\begin{enumerate}

\item[(i)] {\it Let $P \geq 5$ (alt. $P \geq 3$) be any odd prime and let $m
\geq 2$ (alt. $m \geq 3$) 
be any positive integer. Let $q =P^{m}$. Then there exist at least two
$(q+W-6)$-dimensional orbifold lens spaces with fundamental groups of 
order $P^{m}$ which are isospectral but not isometric.}

\item[(ii)] {\it Let $P_{1},P_{2}$ be two odd primes such that $q=P_{1}\cdot
P_{2} \geq 33$. Then there exist at least two $(q+W-6)$-dimensional orbifold 
lens spaces with fundamental groups of order $P_{1} \cdot P_{2}$ which are
isospectral but not isometric.} 

\item[(iii)] {\it Let $q = 2^{m}$ where $m \geq 6$ is any positive integer. Then
there exist at least two $(q+W-5)$-dimensional orbifold lens spaces with
fundamental groups of order $2^{m}$ which are isospectral but not isometric. }

\item[(iv)] {\it Let $q=2P$, where $P \geq 7$ is an odd prime. Then there exist
at least two $(q+W-5)$-dimensional orbifold lens spaces with 
fundamental groups of order $2P$ which are isospectral but not isometric.} 

\end{enumerate}

A consequence of this theorem is that for every integer $x \geq 9$, we can find
a pair of isospectral non-isometric orbifold lens spaces of dimension $x$. 

In Section 4 we look at specific examples that show what the different
generating functions would look like and the types of orbifold lens 
spaces that correspond to each generating function. 

\subsection{Acknowledgments}
I would like to thank Carolyn Gordon very much for some very helpful discussions and 
suggestions, and for reading an early version.


\section{Orbifold Lens Spaces}\label{orblenschapt}

In this section we will generalize the idea of manifold Lens spaces to orbifold
Lens spaces. Manifold Lens spaces are spherical space 
forms where the $n$-dimensional sphere $S^n$ of constant curvature $1$ is acted
upon by a cyclic group of fixed point free isometries 
on $S^n$. We will generalize this notion to orbifolds by allowing the cyclic
group of isometries to have fixed points. For a more general definition of Orbifolds
see Satake \cite{Sat} and Scott \cite{Sc}. For details of spectral geometry on Orbifolds, 
see Stanhope \cite{S1} and E. Dryden, C. Gordon, S. Greenwald and D. Webb in \cite{DGGW}). 

To obtain our main results we will focus on a special subfamily of Lens spaces.
Our technique will parallel Ikeda's technique as 
developed in \cite{I1}. 

\subsection{Preliminaries} \label{subsection31} 

Let $q$ be a positive integer that is not prime. Set \[ q_{0} = \begin{cases}
													\frac{q-1}{2} & \qquad \text{if } q \text{ is odd} \\
													\text{ } \frac{q}{2} & \qquad \text{if } q \text{ is even}.
											    \end{cases} \]
Throughout this article we assume that $q_0 \geq 4$ and that $q$ is not prime. 

For any positive integer $n$ with $2 \leq n \leq q_{0} -2$, we denote by
$\widetilde{I}(q,n)$ the set of $n$-tuples $(p_{1},\ldots,p_{n})$ of 
integers. We define a subset $\widetilde{I}_{0}(q,n)$ of $\widetilde{I}(q,n)$ as
follows: 

\[ \widetilde{I}_{0}(q,n) = \Big \{ (p_1, \ldots, p_n) \in \widetilde{I}(q,n)
\Big \arrowvert p_{i} \not \equiv \pm p_{j} \, (\text{mod} \, q),
1\leq 
i < j \leq n \text{, } g.c.d.(p_1,\dots, p_n, q) =1\Big \}. \]

We introduce an equivalence relation in $\widetilde{I}(q,n)$ as follows:
$(p_1,\ldots,p_n)$ is equivalent to $(s_1,\ldots,s_n)$ if and only if 
there is a number $l$ prime to $q$ and there are numbers $e_i \in \{-1,1\}$ such
that $(p_1,\ldots,p_n)$ is a permutation of $(e_{1}ls_{1},\ldots, e_{n}ls_{n})
\, 
(\text{mod} \, q)$. This equivalence relation also defines an equivalence
relation on $\widetilde{I}_{0}(q,n)$. 

We set $I(q,n) = \widetilde{I}(q,n)/ \sim$ and $I_{0}(q,n) =
\widetilde{I}_{0}(q,n) / \sim$. 

Let $k = q_{0}-n$. We define a map $w$ of $I_{0}(q,n)$ into $I_{0}(q,k)$ as
follows: 

For any element $(p_1,\ldots , p_n) \in \widetilde{I}_{0}(q,n)$, we choose an
element $(q_1,\ldots,q_k) \in \widetilde{I}_{0}(q,k)$ such that 
the set of integers 
	\[ \Big \{ p_1,-p_1, \ldots, p_n,-p_n, q_1,-q_1, \ldots, q_k, -q_k \Big \} \]
forms a complete set of incongruent residues $(\text{mod} \, q)$. Then
we define \[ w( [p_1,\ldots,p_n ]) = [q_1,\ldots,q_k ] \] It is easy to see that 
this map is a well defined bijection.

 The following proposition is similar to a result in \cite{I1}: 
 
 \begin{prop}\label{ch3equivrelation} Let $I_{0}(q,n)$ be as above. Then, \[ |
I_{0}(q,n) | \geq \frac{1}{q_{0}} \binom{q_{0}}{n} \] 
 where $\binom{q_{0}}{n} =1$ if  $q_{0}n = 0$, and \[ \binom{q_{0}}{n} =
\frac{q_{0}!}{n!(q_{0}-n)!} \text{ otherwise.} \]
 \end{prop} 
 
 \begin{proof} 
 	Let $\iqn$ be as above. Consider a subset $\widetilde{I}_{0}'(q,n)$ of $\itqn$
as follows: 
	\[ \widetilde{I}_{0}'(q,n) = \Big \{ (p_1,\ldots, p_n) \in \itqn \Big
\arrowvert \text{at least one of the } p_{i} \text{ is co-prime to } q \Big \}
\] 
	It is easy to see that the equivalence relation on $\itqn$ induces an
equivalence relation on $\widetilde{I}_{0}'(q,n)$. Since we 
	eliminate classes where none of the $p_{i}$'s is co-prime to $q$, we get \[ |
\iqn| \geq | I_{0}'(q,n) | \] where 
	$I'_{0}(q,n) = \widetilde{I}_{0}'(q,n) / \sim$. Now consider a subset $\widetilde{I}_{0}''(q,n)$ of 
	$\widetilde{I}_{0}'(q,n)$ as follows: \[ \widetilde{I}_{0}''(q,n) = \Big \{
(p_1, \ldots, p_n) \in \widetilde{I}_{0}'(q,n) \Big \arrowvert 
	1 = p_1 < \cdots < p_n \leq q_0 \Big \} \] Then it is easy to see that any
element of $\widetilde{I}_{0}'(q,n)$ has an equivalent 
	element in $\widetilde{I}_{0}''(q,n)$. On the other hand, for any equivalence
class in $I_{0}'(q,n)$, the number of elements in 
	$\widetilde{I}_{0}''(q,n)$ which belong to that class is at most $n$. Hence we
have: 
	\[  | \iqn| \geq | I_{0}'(q,n) | \geq \frac{1}{n} \big \arrowvert
\widetilde{I}_{0}''(q,n) \big \arrowvert = \frac{1}{n} \binom{q_{0}-1}{n-1} =
\frac{1}{q_{0}} 	\binom{q_0}{n} \]
	This proves the proposition. 
\end{proof}

\begin{lemma} Let $q = p^m$ or $q = p_{1} \cdot p_2$, where $p,p_1, p_2$ are
primes. Let $D$ be the set of all non-zero integers mod $q$ that 
are not co-prime to $q$. Then $|D|$ is even if $q$ is odd and $|D|$ is odd if
$q$ is even. 
\end{lemma}

\begin{proof} For $q = p^m$.

	If $q$ is odd, then $p$ is an odd prime. $\frac{q}{p} = p^{m-1}$ which is an
odd number. Therefore the number of 
	elements in $D$, $(p^{m-1} - 1)$ is even. 
	
	If $q$ is even, then $p=2$. $\frac{q}{p} = 2^{m-1}$ is even. So the number of
elements in $D$, $(2^{m-1} -1)$, is odd. \\

	For $q = p_{1} \cdot p_{2}$. $(p_{1} \neq p_{2})$
	
	If $q$ is odd, then both $p_1$ and $p_2$ are odd primes. The number of elements
in $D$ is $(\frac{q}{p_1} + \frac{q}{p_2} - 2) = 
	(p_2 + p_1 -2)$ which is even since $p_1+p_2$ is even. 
	
	If $q$ is even, then one of the $p_i$'s is $2$ and the other is an odd prime.
Assume $p_1 = 2$. So, the number of elements in $D$ is 
	$(\frac{q}{p_1} + \frac{q}{p_2} - 2) = (p_2 + p_1 -2) = (p_2 + 2 - 2) = p_2$,
which is odd. 
	
	This proves the lemma.
\end{proof}

We will say that $|D| = 2r$ if $|D|$ is even; and $|D| = 2r-1$ if $|D|$ is odd,
where $r$ is some positive integer. It is easy to see that if $|D|$ is 
even, then exactly $r$ members of $D$ are less than $q_0$. If $|D|$ is odd, 
then $r-1$ members of $D$ are strictly less than $q_0$ and one member of $D$ is
equal to $q_0$ (recall that for even $q$, we set $q_0 = q / 2$, 
and for odd $q$, we set $q_0 = (q-1) / 2$).

With these results we now obtain a better lower bound for $| \iqn |$. 

\begin{prop} \label{prop313} 
	Let $\iqn , I_{0}'(q,n), \widetilde{I}_{0}'(q,n)$ and
$\widetilde{I}_{0}''(q,n)$ be as in \propref{ch3equivrelation}. Let $k = q_0 -
n$. Then 
	\[ | \iqn | \geq \sum_{t=u}^{r} \frac{1}{n-t} \binom{q_{0}-1-r}{n-1-t}
\binom{r}{t} \]
	where $u = r - k$ if $r > k$ and $u =0$ if $r \leq k$, and $r$ is as defined
above. 
\end{prop}

\begin{proof}

The number of ways in which we can assign values to the $p_{i}$'s in $(1=p_1,
p_2, \ldots, p_n) \in \widetilde{I}_{0}''(q,n)$ such that $t$ of the 
$p_i$'s are \emph{not} co-prime to $q$ is \[ \binom{q_{0}-1-r}{n-1-t}
\binom{r}{t} \] 

On the other hand for any equivalence class in $I_{0}'(q,n)$ 
with $t$ of the $p_i$'s not being co-prime to $q$, the number of elements which
belong to that class is at most $n-t$. So the number of such 
possible classes is at least \[ \frac{1}{n-t} \binom{q_{0}-1-r}{n-1-t}
\binom{r}{t} \] Now if $r > k$, this would mean that $n > q_0 - r$, or 
$n -1 > q_0 - 1 -r$. This means that $t$ cannot take any values less than $r-k$,
since that would mean that we are choosing $(n-1-t)$, a 
number larger than $(q_0 -1 -r)$ from $q_0 -1 -r$ and that is not possible. So,
the smallest value for $t$ in this case can be $r - k$. 

On the other hand, if $r \leq k$, then $n \leq q_0 - r$, or $n-1 \leq q_0 -1
-r$. This means that it is possible for us to choose $n$-tuples in 
$\widetilde{I}_{0}''(q,n)$ with all values being co-prime to $q$. Thus, the
smallest value for $t$ would be $0$ in this case. 

It is obvious that the maximum value $t$ can take is $r$ since $(1,p_2, \ldots,
p_n)$ cannot have more than $r$ values that are not co-prime to 
$q$. Now, adding up all the degrees for different values of $t$ we get 
\[ | \iqn | \geq | I_{0}'(q.n) | \geq \sum_{t=u}^{r} \frac{1}{n-t} 
\binom{q_{0}-1-r}{n-1-t} \binom{r}{t} \]
where $u = 0$ if $r \leq k$ and $ u = r-k$ if $r > k$. 

This proves the proposition. 
\end{proof}

\begin{dfn}\label{bj}
\begin{enumerate}

\item [(i)] Let $q$ be a positive integer and $\gamma$ a primitive $q$-th root of $1$.
We denote by $\mathbf{Q}(\gamma)$ the $q$-th cyclotomic field over the rational
number field $\mathbf{Q}$ and denote by $\Phi_{q}(z)$ the $q$-th cyclotomic
polynomial \[ \Phi_{q}(z) = \sum_{t=0}^{q-1} z^{t} \] 
Let $A$ be the set of residues mod $q$ that are co-prime to $q$. We define a map
$\psi_{q,k}$ of $I_{0}(q,k)$ into $\mathbf{Q}(\gamma)[z]$ as 
follows:

 For any equivalence class in $I_{0}(q,k)$, we take an element $(q_1, \ldots,
q_k)$ of $\widetilde{I}_{0}(q,k)$ which belongs to that class. We define
\[ \psi_{q,k}( [q_1, \ldots, q_k ] )(z) = \sum_{l \in A} \, \prod_{i=1}^{k} (z -
\gamma^{q_{i} l})(z - \gamma^{-q_{i}l}) \] This polynomial in 
$\mathbf{Q}(\gamma)[z]$ is independent of the choice of elements which belong to
the class $[q_1, \ldots, q_k ]$. Therefore, the map is well-defined. 

\item [(ii)] Given $q = p^m$, we define 

\[B_j = \big \{ x(mod q)\in\mathbf{Z}^+: p^j\mid x, p^{j+1}\nmid x \big \}\]\
 
We define the maps $\alpha_{q,k}^{(j)}$ of $I_{0}(q,k)$ into
$\mathbf{Q}(\gamma)[z]$ as follows: 

\[ \alpha_{q,k}^{(j)} ( [q_1, \ldots, q_k ] )(z) = \sum_{l \in B_{j}}
\, \prod_{i=1}^{k} (z - \gamma^{q_{i} l})(z - \gamma^{-q_{i}l}) \]

\item [(iii)] Now assume $q = p_{1} \cdot p_{2}$. We define the following sets of
numbers that are not co-prime to $q$. 
\[B = \big \{ xp_{1} \big \arrowvert x = 1, 2, \ldots, (p_{2}-1) \big \} \text{
and } C = \big \{ xp_{2} \big \arrowvert x = 1, 2, \ldots, (p_{1}-1) \big \} \]
We define maps $\alpha_{q,k}$ and $\beta_{q,k}$ as follows:  
\[ \alpha_{q,k} ([q_1. \ldots, q_k] )(z) = \sum_{l \in B} \, \prod_{i=1}^{k} (z
- \gamma^{q_{i} l}) (z - \gamma^{-q_{i} l}) \]
and, 
\[ \beta_{q,k} ([q_1. \ldots, q_k] )(z) = \sum_{l \in C} \, \prod_{i=1}^{k} (z -
\gamma^{q_{i} l}) (z - \gamma^{-q_{i} l}) \]
 
 \end{enumerate}
 \end{dfn}
 
 Since $(z - \gamma^{q_{i} l}) (z - \gamma^{-q_{i} l}) =
(\gamma^{q_{i}l}z-1)(\gamma^{-q_{i}l}z-1)$, the following proposition is easy to
see. 
 
\begin{prop}\label{prop314} If we put 
\begin{align*} 
  	\psi_{q,k} ([q_1, \ldots, q_k])(z) &= \sum_{i=0}^{2k} (-1)^{i} \, a_{i} z^{2k-i}
\\
	\alpha_{q,k}^{(j)}([q_1, \ldots, q_k])(z) &= \sum_{i=0}^{2k} (-1)^{i} \, b_{i,j}
z^{2k-i} \\
	\alpha_{q,k} ([q_1, \ldots, q_k])(z) &= \sum_{i=0}^{2k} (-1)^{i} \, b_{i} z^{2k-i}
\\
	\beta_{q,k} ([q_1, \ldots, q_k])(z) &= \sum_{i=0}^{2k} (-1)^{i} \, c_{i} z^{2k-i} 
\end{align*}
 then we have 
 \begin{enumerate}
 	\item [(i)] $a_i = a_{2k-i}$, $b_{i,j} = b_{(2k-i),j}$, $b_i = b_{2k-i}$ and $c_i =
c_{2k-i}$
	\item [(ii)] $a_0 = |A|$, $b_{0,j} = | B_{j} |$, $b_0 = |B|$ and $c_0 = |C|$. 
\end{enumerate}
\end{prop}

\subsection{Orbifold Lens Spaces and their Generating Functions}\label{subsection32}

	Let $q$ be a positive integer and $p_1,\ldots , p_n$ be $n$ integers mod $q$ such that $g.c.d.(p_1,\dots, p_n, q) =1$. We denote by $g$ the orthogonal matrix given by 
\[ g = 
\begin{pmatrix}
 R(p_{1} / q) & & \text{ {\huge 0}} \\
  & \ddots & \\
  \text{ {\huge 0}} & & R(p_{n} / q)
\end{pmatrix}
\]
where $R(\theta ) = \begin{pmatrix} \cos 2\pi \theta & \sin 2\pi \theta \\ -
\sin 2\pi \theta & \cos 2\pi \theta \end{pmatrix}$
Then $g$ generates the cyclic subgroup $G = \big \{ g^{l} \big \}_{l=1}^{q}$ of
order $q$ of the orthogonal group $O(2n)$. 

We define a Lens space $L(q : p_{1}, \ldots, p_n)$ as follows: 
\[ L(q : p_{1}, \ldots, p_n) = S^{2n-1} / G \] 

$L(q : p_{1}, \ldots, p_n)$ is a good smooth orbifold with $S^{2n-1}$ as its
covering manifold. Let $\pi$ be the covering projection of $S^{2n-1}$ 
onto $S^{2n-1} / G$ 
\[ \pi : S^{2n-1} \rightarrow S^{2n-1} / G \] 
Since the round metric of constant curvature one on $S^{2n-1}$ is
$G$-invariant, it induces a Riemannian metric on $S^{2n-1}/G$. 

We emphasize that, in contrast to the classical notion of Lens space, we do not require 
that the integers $p_i$ be coprime to $q$, and thus we allow the Lens spaces to have 
singular points.  In particular they are good orbifolds.  Henceforth, the term 
"Lens space" will refer to this generalized definition.

\begin{prop} \label{lensspaceisometric}
Let $L = L(q : p_{1}, \ldots, p_n)$ and $L' = L(q : s'_{1}, \ldots, s'_n)$ be
Lens spaces. Then $L$ is isometric to $L'$ if and only if there is a number 
$l$ coprime with $q$ and there are numbers $e_{i} \in \{-1,1 \}$ such that $(p_1, \ldots, p_n)$
is a permutation of $(e_{1}ls_1, \ldots, e_{n}ls_n) \pmod{q}$. 
\end{prop}

\begin{proof} 
Assume $L$ is isometric to $L'$. Let $\phi$ be an isometry between $L$ and $L'$.
Then $\phi$ lifts to an isometry $\tilde{\phi}$ of $S^{2n-1}$. 
In other words, $\phi$ is an orthogonal transformation that conjugates $G$ and
$G'$ where $L = S^{2n-1}/G$ and $L' = S^{2n-1}/G'$. 

So $\tilde{\phi}$ takes $g$, a generator of $G$, to $g'^l$, a generator of $G'$.
This means that the eigenvalues of $g$ and $g'^l$ are the same. This means that
each $p_i$ is equivalent to some $ls_j$ or $-ls_j \pmod{q}$. That means $(p_1,
\ldots, p_n)$ is a permutation of 
 $(e_{1}ls_1, \ldots, e_{n}ls_n) \pmod{q}$, for $e_i \in \{-1,1\}$, $(i = 1,
\ldots, n)$. 
 
 Conversely, assume that there exists an integer $l$ coprime with $q$ and numbers $e_i \in
\{-1,1\}$ $(i = 1, \ldots, n)$ such that $(p_1, \ldots, p_n)$ is a 
 permutation of $(e_{1}ls_1, \ldots, e_{n}ls_n) \pmod{q}$. 
 
 Note that the isometry of $S^{2n-1}$ onto $S^{2n-1}$ defined by the map 
 \[ (z_1, \ldots, z_i, \ldots, z_n) \longmapsto (z_{ \sigma (0) }, \ldots ,
\bar{z}_{ \sigma (i) }, \ldots , z_{ \sigma (n) } )    \]
 where $\sigma$ is a permutation, induces an isometry of $L(q: s_1, \ldots,
s_n)$ onto \\
 $L(q: s_{\sigma (1)} , \ldots , s_{\sigma (i)}, \ldots,  s_{\sigma (n)})$.
Since $g'^l$ is also a generator of $G'$, the Lens space 
 $L(q: ls_1, \ldots, ls_n)$ is identical to $L(q: s_1, \ldots, s_n)$. 
 
 Now the above isometry induces an isometry of $L(q: s_1, \ldots, s_n)$ onto
$L(q: e_{1}ls_1, \ldots, e_{n}ls_{n})$. This means that 
 $L'  = L(q: s_1, \ldots, s_n)$ is isometric to $L(q: e_{1}ls_1, \ldots,
e_{n}ls_{n})$. But $(e_{1}ls_1, \ldots, e_{n}ls_{n})$ is simply a 
 permutation of $(p_1,\ldots, p_n)\pmod{q}$. Therefore $L(q: p_1, \ldots, p_n)$
is isometric to $L(q: e_{1}ls_1, \ldots, e_{n}ls_{n})$, which, 
 in turn, is isometric to $L(q: s_1, \ldots, s_n)$. This proves the converse. 
 \end{proof}
 
 For any $f \in C^{\infty} (\lqp)$, we define the Lapacian on the lense space as
$\widetilde{\Delta}(\pi^{*} f) = \pi^{*} (\Delta f)$. 
We construct the generating function associated with the Laplacian on $\lqp$ analogous to 
 the construction in the manifold case (see \cite{I1}, \cite{I2} and \cite{IY}). 
 
 Let $\tilde{\Delta}$, $\Delta$ and $\Delta_0$ denote the Laplacian of
$S^{2n-1}$, $\lqp$ and $\mathbf{R}^{2n}$, respectively. 
 
 \begin{dfn} For any non-negative real number $\lambda$, we define the
\emph{eigenspaces} $\widetilde{E}_{\lambda}$ and $E_{\lambda}$ as follows: 
 \begin{align*}
 	\widetilde{E}_{\lambda} &= \big \{ f \in C^{\infty} (S^{2n-1}) \big \arrowvert
\widetilde{\Delta} f = \lambda f \big \} \\ 
	E_{\lambda} &= \big \{ f \in C^{\infty} (\lqp) \big \arrowvert \Delta f =
\lambda f \big \}
\end{align*} 
\end{dfn} 

The following lemma follows from the definitions of $\Delta$ and smooth function. 

\begin{lemma}
Let $\big (\widetilde{E}_{\lambda} \big )_{G}$ be the space of all $G$-invariant
functions of $\widetilde{E}_{\lambda}$. Then $\dim (E_{\lambda}) = 
\dim (\widetilde{E})_{G}$. 
\end{lemma} 

Let $\Delta_0$ be the Laplacian on $\mathbf{R}^{2n}$ with respect to the flat
K\"ahler metric. Set $r^2 = \sum_{i=1}^{2n} x_{i}^2$, where $(x_1, x_2, \ldots, 
x_{2n})$ is the standard coordinate system on $\mathbf{R}^{2n}$. 

For $k \geq 0$, let $P^k$ denote the space of complex valued homogenous
polynomials of degree $k$ on  $\mathbf{R}^{2n}$. Let $H^k$ be the subspace 
of $P^k$ consisting of harmonic polynomials on $\mathbf{R}^{2n}$, 
\[ H^k = \big \{ f \in P^k \big \arrowvert \Delta_{0}f = 0 \big \} \]
Each orthogonal transformation of $\mathbf{R}^{2n}$ induces canonically a linear
isomorphism of $P^k$.  

\begin{prop} \label{ch3inducedlinmap}
The space $H^k$ is $O(2n)$-invariant, and $P^k$ has the direct sum
decomposition: $P^k = H^k \oplus r^{2} P^{k-2}$. 
The injection map $i: S^{2n-1} \rightarrow \mathbf{R}^{2n}$ induces a linear map
$i^{*}: C^{\infty}(\mathbf{R}^{2n}) \rightarrow C^{\infty}(S^{2n-1})$. 
We denote $i^{*}(H^k)$ by $\mathcal{H}^k$.
\end{prop}

\begin{prop} \label{ch3eigenspaceprop}
$\mathcal{H}^k$ is an eigenspace of $\widetilde{\Delta}$ on $S^{2n-1}$ with
eigenvalue $k(k+2n-2)$ and $\sum_{k=0}^{\infty} \mathcal{H}^{k}$ 
is dense in $C^{\infty}(S^{2n-1})$ in the uniform convergence topology.
Moreover, $\mathcal{H}^k$ is isomorphic to $H^k$. That is, 
$i^{*}\! : H^k \stackrel{\simeq}{\longrightarrow} \mathcal{H}^k$. 
\end{prop}

For proofs of these propositions, see \cite{BGM}. 

Now by \corref{ch3dimcorollary} and \propref{ch3eigenspaceprop}, we get

\begin{cor}\label{corollary327} Let \[ L(q : p_{1}, \ldots, p_n) = S^{2n-1} / G \] be a lens space and
$\mathcal{H}_{G}^k$ be the space of all $G$-invariant functions in
$\mathcal{H}^k$. Then 
\[ \dim E_{k(k+2n-2)} = \dim \mathcal{H}_{G}^k \]
Moreover, for any integer $k$ such that $\dim \mathcal{H}_{G}^k \neq 0$,
$k(k+2n-2)$ is an eigenvalue of $\Delta$ on $\lqp$ with multiplicity equal to 
$\dim \mathcal{H}_{G}^k$, and no other eigenvalues appear in the spectrum of
$\Delta$. 
\end{cor}

\begin{dfn}
The \emph{generating function} $F_{q}(z : p_1, \ldots, p_n)$ associated to the
spectrum of the Laplacian on $\lqp$ is the generating function 
associated to the infinite sequence $\big \{ \dim \mathcal{H}_{G}^{k} \big
\}_{k=0}^{\infty}$ , i.e., 
\[ F_{q}(z : p_1, \ldots, p_n) = \sum_{k =0}^{\infty} \big( \dim
\mathcal{H}_{G}^{k} \big) z^k \] 
\end{dfn}

By \corref{ch3dimcorollary} we know that the generating function determines the
spectrum of $\lqp$. This fact gives us the following proposition:

\begin{prop} \label{prop329}
Let $L = \lqp$ and $L' = L(q': s_1, \ldots, s_n)$ be two lens spaces. Let
$F_{q}(z: p_1, \ldots, p_n)$ and $F_{q'}(z: s_1, \ldots, s_n)$ be the generating
functions
associated to the spectrum of $L$ and $L'$, respectively. Then $L$ is
isospectral to $L'$ if and only if $F_{q}(z: p_1, \ldots, p_n) = F_{q'}(z: s_1,
\ldots, s_n)$. 
\end{prop}

The following theorem (\cite{I1} and \cite{I2}) holds true for the orbifold lens spaces as well.

\begin{thm} \label{theorem3210}
Let $\lqp$ be a lens space and $F_{q}(z: p_1, \ldots, p_n)$ the generating
function associated to the spectrum of $\lqp$. Then, on the domain 
$\big \{ z \in \mathbb{C} \, \big \arrowvert  \, |z| < 1 \big \}$, 
\[ F_{q}(z: p_1, \ldots, p_n) = \frac{1}{q} \, \sum_{l=1}^{q} \,
\frac{1-z^2}{\prod_{i=1}^{n} (z - \gamma^{p_{i} l} )(z - \gamma^{-p_{i} l})} \] 
\end{thm}

\begin{cor}
Let $\lqp$ be isospectral to $L(q' : s_1, \ldots, s_n)$. Then $q = q'$. 
\end{cor}

Now let $\widetilde{\mathcal{L}} (q,n) $ be the family of all
$(2n-1)$-dimensional lens spaces with fundamental groups of order $q$, and let 
$\widetilde{\mathcal{L}}_{0} (q,n) $ be the subfamily of
$\widetilde{\mathcal{L}} (q,n) $ defined by: 
\[ \widetilde{\mathcal{L}}_{0} (q,n) = \big \{ L(q: p_1, \ldots, p_n) \in
\widetilde{\mathcal{L}} (q,n) \big \arrowvert p_i \not \equiv \pm p_j \! \! \!
\pmod{q}, 
	1\leq i < j \leq n \big \} \] 
The set of isometry classes of $\widetilde{\mathcal{L}} (q,n) $ is denoted by $
\mathcal{L} (q,n) $, and the set of isometry classes of 
$ \widetilde{\mathcal{L}}_{0} (q,n) $ is denoted by $\mathcal{L}_{0} (q,n) $. 

By \propref{lensspaceisometric}, the map $\lqp \mapsto (p_1, \ldots, p_n)$ of
$\widetilde{\mathcal{L}}_{0} (q,n)$  $\big [ \text{resp. }
\widetilde{\mathcal{L}}(q,n) \big ]$ onto $\widetilde{I}_{0}(q,n)$ $ \big [
\text{resp. } \widetilde{I} (q,n) \big ]$ induces a one-to-one map between
$\mathcal{L}_{0}(q,n)$ and $I_{0}(q,n)$ 
$\big [ \text{resp. } \mathcal{L}(q,n) \text{ and } I(q,n) \big ]$.  

The above fact, together with \propref{prop313}, gives us the following: 

\begin{prop} \label{prop3212}
Retaining the notations as above, we get 
\[ \abs{\mathcal{L}_{0}(q,n)} \geq \sum_{t=u}^{r} \frac{1}{n-t}
\binom{q_{0}-1-r}{n-1-t} \binom{r}{t}   \] 
where $u = r- k$ if $r > k$, and $u = 0$ if $r \leq k$; $r$ is the number of
residues $\mod{q}$ that are not co-prime to $q$ and are less than or equal to
$q_0$. 
\end{prop} 

Note that by \propref{ch3equivrelation}, we also get that 
\[ \abs{\mathcal{L}_{0}(q,n)} \geq \frac{1}{q_0} \binom{q_0}{n} \] 
Next, we will re-formulate the generating function $F_{q}(z: p_1, \ldots, p_n)$
in a form that will help us find isospectral pairs that are non-isometric. 

\begin{prop} \label{prop3213}
Let $\lqp$ be a lens space belonging to $\widetilde{\mathcal{L}}_{0}(q,n)$, $k =
q_0 - n$, and let $w$ be the map of $I_0 (q,n)$ onto $I_0 (q,k)$ defined 
in section \ref{subsection31}. Then 
\begin{enumerate}
	\item[(i)] If $q = P^m$, where $P$ is a prime, we have 
		\begin{multline*}
		F_{q}(z: p_1, \ldots, p_n) = \frac{1}{q} \Bigg \{ \frac{ (1-z^2) } {
(1-z)^{2n}} + \frac{ \psi_{q,k} (w([p_1, \ldots, p_n]))(z) (1-z^2) } { \Phi_{q} (z) }
\, + \\
		\sum_{j=1}^{m-1} \, \frac{ \alpha_{q,k}^{(j)} (w([p_1, \ldots, p_n]))(z) (1-z^2)
}{ ( \Phi_{P^{m-j}} (z) )^{P^j} (1 - z)^{P^j - 1}} \Bigg \} 
		\end{multline*} 
	
	\item[(ii)] If $q = P_1 \cdot P_2 $, where $P_1$ and $P_2$ are primes, we have 
		\begin{multline*}
		F_{q}(z: p_1, \ldots, p_n) = \frac{1}{q} \Bigg \{ \frac{ (1-z^2) } {
(1-z)^{2n} } + \frac{ \psi_{q,k} (w([p_1, \ldots, p_n]))(z) (1-z^2) } { \Phi_{q} (z)
}  \\ 
		+ \, \frac{ \alpha_{q,k}  (w( [p_1, \ldots, p_n] ))(z) (1-z^2 ) } { ( \Phi_{P_2}
(z) )^{P_{1}} (1 - z)^{P_{1} - 1} }  \\
		+ \, \frac{ \beta_{q,k}  (w( [p_1, \ldots, p_n] ))(z) (1-z^2 ) } { ( \Phi_{P_1}
(z) )^{P_{2}} (1 - z)^{P_{2} - 1} }  \Bigg \} 
		\end{multline*} 	
\end{enumerate} 
where $\psi_{q,k}$, $\alpha^{(j)}_{q,k}$, $\alpha_{q,k}$ and $\beta_{q,k}$ are
as defined in definition \ref{bj} and $\Phi_{t}(z) = \sum_{v = 0}^{t-1} z^v$.  
\end{prop} 

\begin{proof}
We choose integers $q_1, \ldots, q_k$ such that the set of integers 
\[ \{ p_1, -p_1, \ldots, p_n, -p_n, q_1, -q_1, \ldots, q_k, -q_k \} \]
forms a complete set of residues mod $q$. 
\begin{enumerate}
	\item[(i)] We write 
	\begin{multline*}
	F_{q}(z: p_1, \ldots, p_n) = \frac{1}{q} \Bigg [ \sum_{l \in A} \frac{ (1 -
z^2) }{ \prod_{i = 1}^{n} (z - \gamma^{p_{i}l})(z - \gamma^{-p_{i}l})} \\
	+ \, \sum_{j=1}^{m-1} \sum_{l \in B_{j}} \frac{ (1-z^2) } { \prod_{i = 1}^{n}
(z - \gamma^{p_{i}l})(z - \gamma^{-p_{i}l})} \Bigg \}
	\end{multline*}
	Now, for any $l \in A$, we have
	\begin{align*}
	\frac{1}{ \prod_{i = 1}^{n} (z - \gamma^{p_{i}l})(z - \gamma^{-p_{i}l})} =
\frac{ \prod_{i=1}^{k} (z - \gamma^{q_{i}l})(z - \gamma^{-q_{i}l})}{\Phi_{q}(z)}
	\end{align*}
	For $l \in B_{j}$, we have 
	\begin{align*}
	\frac{1}{ \prod_{i = 1}^{n} (z - \gamma^{p_{i}l})(z - \gamma^{-p_{i}l})} =
\frac{ \prod_{i = 1}^{k}  (z - \gamma^{q_{i}l})(z - \gamma^{-q_{i}l})}{
	(\Phi_{P^{m-j}}(z))^{P^j} (1 - z)^{P^j -1}}
	\end{align*}
	where $\Phi_{t}(z) = \sum_{v = 0}^{t-1} z^v$.
	 
	Now, $(i)$ follows from these facts.
	\item[(ii)] We write
	\begin{multline*}
	F_{q}(z: p_1, \ldots, p_n) =  \frac{1}{q} \Bigg [ \sum_{l \in A} \frac{ (1 -
z^2) }{ \prod_{i = 1}^{n} (z - \gamma^{p_{i}l})(z - \gamma^{-p_{i}l})} \\ 
	+ \, \sum_{l \in B} \frac{ (1 - z^2) }{ \prod_{i = 1}^{n} (z -
\gamma^{p_{i}l})(z - \gamma^{-p_{i}l})}
         + \, \sum_{l \in C} \frac{ (1 - z^2) }{ \prod_{i = 1}^{n} (z -
\gamma^{p_{i}l})(z - \gamma^{-p_{i}l})} \Bigg ]
         \end{multline*}
         Again, for $l \in A$, 
         \begin{align*} 
         \frac{1}{ \prod_{i = 1}^{n} (z - \gamma^{p_{i}l})(z -
\gamma^{-p_{i}l})} = \frac{ \prod_{i=1}^{k} (z - \gamma^{q_{i}l})(z -
\gamma^{-q_{i}l})}{\Phi_{q}(z)}
         \end{align*}
         For $l \in B$, we have
         \begin{align*}
	\frac{1}{ \prod_{i = 1}^{n} (z - \gamma^{p_{i}l})(z - \gamma^{-p_{i}l})} =
\frac{ \prod_{i = 1}^{k}  (z - \gamma^{q_{i}l})(z - \gamma^{-q_{i}l})}{
	(\Phi_{P_2}(z))^{P_1} (1 - z)^{P_1 -1}}
	\end{align*}
	For $l \in C$, we have 
	\begin{align*}
	\frac{1}{ \prod_{i = 1}^{n} (z - \gamma^{p_{i}l})(z - \gamma^{-p_{i}l})} =
\frac{ \prod_{i = 1}^{k}  (z - \gamma^{q_{i}l})(z - \gamma^{-q_{i}l})}{
	(\Phi_{P_1}(z))^{P_2} (1 - z)^{P_2 -1}}	
         \end{align*}
         Now, $(ii)$ follows from these facts. 
\end{enumerate}	
\end{proof} 

From \propref{prop329} and \propref{prop3213}, we get the following proposition 
\begin{prop} \label{prop3214}
Let $L = \lqp$ and $L' = L(q: s_1,\ldots, s_n)$ be lens spaces belonging to
$\mathcal{L}_{0}(q,n)$. Let $k = q_0 - n$. 
\begin{enumerate}
\item[(i)] If $q = P^m$, then $L$ is isospectral to $L'$ if 
	\begin{align*}
	& \psi_{q,k}(w([p_1, \ldots, p_n])) = \psi_{q,k}(w([s_1, \ldots, s_n])) \\ 
	\text{and} \qquad & \alpha_{q,k}^{(j)}(w([p_1, \ldots, p_n])) =
\alpha_{q,k}^{(j)}(w([s_1, \ldots, s_n]))
	\end{align*}
	for $j =1, \ldots, m-1$. 
	
\item[(ii)] If $q = P_1 \cdot P_2$, then $L$ is isospectral to $L'$ if 
	\begin{align*}
	& \psi_{q,k}(w([p_1, \ldots, p_n])) = \psi_{q,k}(w([s_1, \ldots, s_n])) \\ 
	& \alpha_{q,k} (w([p_1, \ldots, p_n])) = \alpha_{q,k} (w([s_1, \ldots, s_n]))
\\ 
	\text{and} \qquad & \beta_{q,k} (w([p_1, \ldots, p_n])) = \beta_{q,k} (w([s_1,
\ldots, s_n]))
	\end{align*}
\end{enumerate}
\end{prop}






\section{Isospectral Non-isometric Lens Spaces}

By applying \propref{prop3212} and \propref{prop3214} we will obtain our main
theorem in this Section for odd-dimensional
lens spaces. Next, we will extend the results to obtain even-dimensional pairs
of lens spaces corresponding to every pair 
of odd-dimensional lens spaces. 

\subsection{Odd-Dimensional Lens Spaces}\label{section41}

From the results in Sections 2 and 3 we get the following diagrams: \\
For $q = P^m$, 
\begin{align}
	\mathcal{L}_{0}(q,n) \xrightarrow{\sim} I_{0}(q,n) \xrightarrow[w]{\sim}
I_{0}(q,k) \xrightarrow[ \tau^{(m)}_{q,k}] {} 
	Q^m (\gamma)[z] \label{eq41}
\end{align}
where $\tau^{(m)}_{q,k} = ( \psi_{q,k}, \alpha_{q,k}^{(1)}, \ldots,
\alpha_{q,k}^{(m-1)} )$, and $Q^m (\gamma)[z]$ denotes $m$ copies of the field of rational polynomials $Q (\gamma)[z]$ \\
For $q = P_1 \cdot P_2$, 
\begin{align}
\mathcal{L}_{0}(q,n) \xrightarrow{\sim} I_{0}(q,n) \xrightarrow[w]{\sim} I_0
(q,k) \xrightarrow[ \mathcal{S}^{(3)}_{q,k} ]{}
Q^{3} (\gamma)[z] \label{eq42}
\end{align}
where $\mathcal{S}^{(3)}_{q,k}  = (\psi_{q,k}, \alpha_{q,k}, \beta_{q,k} )$. 

Now, from \propref{prop3214}, if $\tau^{(m)}_{q,k}$ [resp.
$\mathcal{S}^{(3)}_{q,k} $] is not one-to-one, then we will have 
non-isometric lens spaces having the same generating function. This would give
us our desired results. 

We first calculate the values for the required coefficients of 
$\psi_{q,2}$, $\alpha^{(j)}_{q,2}$, 
$\alpha_{q,2}$ and $\beta_{q,2}$. 

Using \propref{prop314} we can calculate the values of the various coefficients of
$\psi_{q,k}$, $\alpha^{(j)}_{q,k}$, 
$\alpha_{q,k}$ and $\beta_{q,k}$.

First we will find coefficients for $z$ and $z^2$ for any given $k$, and from
that we can find the values for $k = 2$. 

From the definitions of $\psi_{q,k} ([q_1, \ldots, q_k])$, in the notation of 2.1.5, it is easy to see that 
	\[ a_1 = \sum_{i=1}^{k} \sum_{l \in A} \gamma^{q_i l} + \sum_{i =1}^{k} \sum_{l
\in A} \gamma^{-q_i l} = 
	2 \sum_{i=1}^{k} \sum_{l \in A} \gamma^{q_i l} \] 
Similarly, 
	\begin{align*}
		b_{1,j} &= 2 \sum_{i=1}^{k} \sum_{l \in B_j} \gamma^{q_i l} \\
		b_1 &= 2 \sum_{i=1}^{k} \sum_{l \in B} \gamma^{q_i l} \\
		c_1 &= 2 \sum_{i=1}^{k} \sum_{l \in C} \gamma^{q_i l} 
	\end{align*}
Also, 
	\begin{align*}
	a_2 &= \sum_{l \in A} \Big [ k+ \sum_{1 \leq i < t \leq k} \gamma^{(q_i + q_t)
l} + \sum_{1 \leq i < t \leq k}
	\gamma^{-(q_i + q_t)l}  + \sum_{1 \leq i < t \leq k} \gamma^{(q_i - q_t) l} +
\sum_{1 \leq i < t \leq k} \gamma^{-(q_i - q_t) l} 
	\Big ] \\
	        &= k \abs{A} + 2 \sum_{l \in A} \sum_{1 \leq i < t \leq k} \gamma^{(q_i
+ q_t) l} + 
	        		 2\sum_{l \in A} \sum_{1 \leq i < t \leq k} \gamma^{(q_i - q_t) l}
	\end{align*}
Similarly, 
	\begin{align*}
	b_{2,j} &= k \abs{B_{j}} + 2 \sum_{l \in B_j} \sum_{1 \leq i < t \leq k}
\gamma^{(q_i + q_t) l} + 
		2 \sum_{l \in B_j} \sum_{1 \leq i < t \leq k} \gamma^{(q_i - q_t) l} \\
	b_2 &=  k \abs{B} + 2 \sum_{l \in B} \sum_{1 \leq i < t \leq k} \gamma^{(q_i +
q_t) l} + 
		2 \sum_{l \in B} \sum_{1 \leq i < t \leq k} \gamma^{(q_i - q_t) l} \\
	c_2 &= k \abs{C} + 2 \sum_{l \in C} \sum_{1 \leq i < t \leq k} \gamma^{(q_i +
q_t) l} + 
		2 \sum_{l \in C} \sum_{1 \leq i < t \leq k} \gamma^{(q_i - q_t) l}
	\end{align*} 
where $\abs{A}$, $\abs{B_j}$, $\abs{B}$ and $\abs{C}$ are cardinalities of $A$,
$B_j$, $B$ and $C$ - as defined in 2.1.4 - respectively. 

\begin{prop} \label{prop411a}
	Let $p$ be an odd prime and let $q = p^m$ where $m$ is an integer greater than
$1$. Let $q_0 = \frac{q-1}{2}$. Let $k=2$ and 
	$n = q_0 -2$. Then the maps $\tau_{q,k}^{(m)}$ and $\mathcal{S}_{q,k}^{(3)}$ as
defined in (\ref{eq41}) and (\ref{eq42}) (and hence the generating function) are 
dependent only on where the various $q_i$'s and their sums and differences reside. 
\end{prop}

In a similar fashion we can find values of coefficients of higher powers of $z$
when $k > 2$. These coefficients will contain 
terms that include higher sums and differences of the various $q_i$'s in the
powers of $\gamma$. 

We will prove two propositions that will give us upper bounds on the
number of expressions for $\tau_{q,k}^{(j)}$ and 
$\mathcal{S}_{q,k}^{(3)}$, respectively, where $k = 2$. 

\begin{prop} \label{prop411}
	Let $p$ be an odd prime and let $q = p^m$ where $m$ is an integer greater than
$1$. Let $q_0 = \frac{q-1}{2}$. Let $k=2$ and 
	$n = q_0 -2$. Then the number of expressions that $\tau^{(j)}_{q,2}$ can have is
at most $m^2$. 
\end{prop}
\begin{proof}
	We will find the number of $\tau^{(j)}_{q,2}$ by considering the following
cases:
	\begin{enumerate}
	\item[\underline{Case 1:}] $q_1, q_2 \in B_j \quad (j = 1, 2, \ldots, (m-1))$
where $B_j = \big \{ x(mod q)\in\mathbf{Z}^+: p^j\mid x, p^{j+1}\nmid x \big \} \\$
	If $q_1,q_2\in B_j$, then either both of $q_1\pm q_2$ lie in $B_j$ or else one
lies in $B_j$ and the other in some $B_k$ with $j<k\leq m-1$.  Thus there are
$m-j$ possibilities.  As $j$ varies from $1$ to $m-1$, we thus obtain 
$(m-1)+(m-2)+\dots +1=\frac{m(m-1)}{2}$
expressions.  

\bigskip
	
	\item[\underline{Case 2:}] $q_1 \in B_j$ and $q_2 \in B_t$, $B_j \neq B_t$ \\
We may assume $j<t$.  It follows that $q_1\pm q_2$ both lie in $B_j$.  Thus as
$j$ and $t$ vary, we obtain ${{m-1} \choose {2}} =\frac{(m-1)(m-2)}{2}$
expressions.

\bigskip
	
	\item[\underline{Case 3:}] $q_1 \in B_j$ and $q_2 \in A$, or vice versa \\
	Here we note that $q_1 \pm q_2$ always belongs to $A$. Therefore, in this case
we will get $(m-1)$ possible expressions for $\tau^{(j)}_{q,2}$,
	one each for the case where $q_1 \in A$ and $q_2 \in B_j \quad (j=1, 2, \ldots,
(m-1))$, or vice versa. 

\bigskip
	
	\item[\underline{Case 4:}] $q_1,q_2 \in A$ \\
	We will get $1$ possible expression if $q_1 \pm q_2 \in A$. Then we will get $1$
possible expression each for the case when $q_1 + q_2 \in 
	A$ and $q_1 - q_2 \in B_j$ (or vice versa) for $j = 1, 2, \ldots, (m-1)$. There
are no other possibilities in this case. So the maximum number 
	of possible expressions for $\tau^{(j)}_{q,2}$ in this case will be $m-1+1 = m$.  
	
\bigskip

	Case $1$ though Case $4$ are the only possible cases that occur for $k=2$.
Adding up the numbers of all possible expressions 
	for $\tau^{(j)}_{q,2}$ from each case we get the maximum number of possible
expressions that $\tau^{(j)}_{q,2}$ can have: 
	\begin{multline*}
	\frac{m(m-1)}{2} + \frac{(m-1)(m-2)}{2} + (m-1) + m \\ 
	= \frac{m^2 - m + m^2 - 3m + 2 + 2m -2 +2m}{2} = \frac{2m^2}{2} = m^2
	\end{multline*}	
	\end{enumerate}	
\end{proof}

\begin{prop} \label{prop412}
	Let $q = p_1 \cdot p_2$, where $p_1, p_2$ are distinct odd primes. Let $q_0 =
\frac{q-1}{2}$. Let $k=2$ and $n = q_0 - 2$. Then the number
	of possible expressions for $\mathcal{S}^{(3)}_{q,2}$ is at most $11$. 
\end{prop}
\begin{proof} 
	As in the previous proposition, we prove this result by considering all the
possible cases for $q_1$ and $q_2$(where $q_1\pm q_2$ is not congruent to $0(mod q)$).  
	\begin{enumerate}
	\item[\underline{Case 1}] $q_1, q_2 \in B \quad (\text{or } q_1, q_2 \in C)$,
where $B = \big \{xp_1 \big \arrowvert x = 1, \ldots,
	(p_2 - 1) \big \}$ and $C = \big \{ xp_2 \big \arrowvert x = 1, \ldots, (p_1 -
1) \big \}$. 
	
	Then $q_1 \pm q_2 \in B$ $(\text{or } q_1 \pm q_2 \in C \text{,
respectively})$. There are no other possibilities for this case. 
	
	\item[\underline{Case 2:}] $q_1 \in B$ and $q_2 \in C$ (or vice versa). \\
	We have just one possible expression in this case, when $q_1 \pm q_2 \in A$. 
	
	\item[\underline{Case 3:}] $q_1 \in A, q_2 \in B$ or $q_1 \in A, q_2 \in C$ (or
vice versa). \\ 
	We will get one expression each when $q_1 \pm q_2 \in A$. Then we will get one
possible expression for the case when $q_1 \in A, q_2 \in B$, 
	and $q_1 + q_2 \in A, q_1 - q_2 \in C$, (or vice versa). \\
	We will get one more possible expression for the case when $q_1 \in A, q_2 \in
C$, and $q_1 + q_2 \in A, q_1 - q_2 \in B$ (or vice versa). \\
	So, in this case we get a possible $4$ expressions for $\mathcal{S}^{(3)}_{q,2}$. 
	
	\item[\underline{Case 4:}] $q_1, q_2 \in A$. \\
	We will get one possible expression where $q_1 \pm q_2 \in A$. 
	We get another possible expression where $q_1 + q_2 \in A$ and $q_1 - q_2 \in B$
(or vice versa). 
	We get a third possible expression where $q_1 + q_2 \in A$ and $q_1 - q_2 \in C$
(or vice versa). 
	We get a fourth possible expression where $q_1 + q_2 \in B$ and $q_1 - q_2 \in C$
(or vice versa). \\
	So, we get a total of $4$ possible expressions for $\mathcal{S}^{(3)}_{q,2}$ in
this case. 
	\end{enumerate}
	
	Case $1$ through Case $4$ are the only possible cases than can occur for $k =
2$. Adding up the number of all possible expressions for 
	$\mathcal{S}^{(3)}_{q,2}$ from each case we get the maximum number of possible
expressions for $\mathcal{S}^{(3)}_{q,2}$: 
	\begin{align*}
	2 + 1 + 4 + 4 = 11
	\end{align*}
\end{proof}

	It is important to note that in the above propositions the number of possible
expressions is the \emph{maximum} number of expressions 
	that can happen. It is possible that for a given $q = p^m$ or $q = p_1 \cdot
p_2$ not all the expressions will occur. We will see this in an
	example later. 
	
	We now prove two similar propositions for even $q$ of the form $2^m$ and $2p$,
where $m$ is a positive integer and $p$ is a 
	prime. 
	
\begin{prop} \label{prop413}
	Let $q = 2^m$ where $m \geq 3$. Let $q_0 = \frac{q}{2}$, i.e., $q_0 = 2^{m-1}$.
Let $k = 2$ and $n = q_0 - 2$. Then the number of 
	possible expressions that $\tau^{(j)}_{q,2}$ can have is at most $(m-1)^2$. 
\end{prop}
\begin{proof}
	We proceed as in the previous propositions. 
	\begin{enumerate}
	\item[\underline{Case 1:}] $q_1 , q_2 \in B_j \quad (j = 1, 2, \ldots, (m-3))$,
where $B_j = \big \{ x(mod q)\in\mathbf{Z}^+: 2^j\mid x, 2^{j+1}\nmid x \big \} \\$
	
	We first note that the cases where $q_1,q_2 \in B_{m-2}$ or $B_{m-1}$ will not
occur.	Now when $q_1, q_2 \in B_j$, then one of the $q_1 + q_2$ or $q_1 - q_2$ will
belong to $B_{j+1}$ and the other will belong to $B_t$
	for $t > j+1$. 
	
	Now, with $q_1, q_2 \in B_j$ (where $j < m-2$), we get $(m-2-j)$ possible expressions for
$\tau^{(j)}_{q,2}$.
	
	So, in this case, the total number of possible expressions for $\tau^{(j)}_{q,2}$
are: 
	\begin{align*}
	(m-3) + (m-4) + \cdots + 3 + 2 + 1 = \frac{(m-2)(m-3)}{2}
	\end{align*}
	
	\item[\underline{Case 2:}] $q_1 \in B_j$ and $q_2 \in B_t$, where $B_j \neq
B_t$. \\
	We can assume that $j < t$. This would mean that $q_1 \pm q_2 \in B_j$ always.
So, as in Case $2$ of \propref{prop411}, we get 
	that the total number of expressions for $\tau^{(j)}_{q,2}$ will be
$\frac{(m-1)(m-2)}{2}$. 
	
	\item[\underline{Case 3:}] $q_1 \in B_j$ and $q_2 \in A$ (or vice versa). \\
	We notice that $q_1 \pm q_2 \in A$ always. So, just like in Case $3$ of
\propref{prop411}, we will get that the total number of possible 
	expressions for $\tau^{(j)}_{q,2}$ will be $(m-1)$. 
	
	\item[\underline{Case 4:}] $q_1,q_2 \in A$. \\
	In this case one of $q_1 + q_2$ or $q_1 - q_2$ will belong to $B_1$ and the
other will belong to one of the $B_j$ for $j > 1$. 
	Therefore, for this case we will get $(m-2)$ possible expressions for
$\tau^{(j)}_{q,2}$, one each for the case when $q_1 + q_2 \in B_1$
	 (alt. $q_1 - q_2 \in B_1$) and $q_1 - q_2 \in B_t$ (alt. $q_1 + q_2 \in B_t$)
for $t = 2,3, \ldots, m-1$. 
	\end{enumerate}
	 
	 Now, adding up all the possible expressions from the four cases above we get the
maximum number of possible expressions for $\tau^{(j)}_{q,2}$: 
	 \begin{align*}
	 &\frac{(m-2)(m-3)}{2} + \frac{(m-1)(m-2)}{2} + (m-1) + (m-2) \\
	 &= \frac{m^2 -5m + 6 + m^2 - 3m + 2 + 2m - 2 + 2m - 4}{2} \\
	 &= m^2 - 2m + 1 = (m-1)^2
	 \end{align*}
	
\end{proof}

Our next proposition gives us the maximum number of expressions for
$\mathcal{S}^{(3)}_{q,2}$ when $q = 2p$ for some
prime $p$. 

\begin{prop} \label{prop414}
	Let $q = 2p$ where $p$ is an odd prime. Let $q_0 = \frac{q}{2} = p$. Let $k =
2$ and $n = q_0 -2$. Then the number of
	possible expressions for $\mathcal{S}^{(3)}_{q,2}$ is at most $6$. 
\end{prop}
\begin{proof} 
	As before we will analyze the different possible cases. Note that in this
situation we have $B = \big \{ 2, 4, 6, 
	\ldots, 2(p-1) \big \}$ and $C = \{p \}$. 
	\begin{enumerate}
	\item[\underline{Case 1:}] $q_1, q_2 \in B$. We will have $q_1 \pm q_2 \in B$
always. \\
	Notice that in this case $q_1,q_2$ cannot belong to $C$ since $C$ has only one
element. So we get $1$ 
	possible expression in this case for $\mathcal{S}^{(3)}_{q,2}$. 
	
	\item[\underline{Case 2:}] $q_1 \in B, q_2 \in C$. In this case $q_1 \pm q_2
\in A$ always. \\
	So, we get $1$ possible expression in this case for $\mathcal{S}^{(3)}_{q,2}$. 
	
	\item[\underline{Case 3:}] $q_1 \in A, q_2 \in B$ or $q_1 \in A, q_2 \in C$. \\
	When $q_1 \in A$ and $q_2 \in C$, then $q_1 \pm q_2 \in B$ always. So, we get
$1$ possible expression for 
	$\mathcal{S}^{(3)}_{q,2}$. When $q_1 \in A, q_2 \in B$, we will get $1$
possible expression for the situation when 
	$q_1 \pm q_2 \in A$. We will get another possible expression for
$\mathcal{S}^{(3)}_{q,2}$ where $q_1 + q_2 \in A$ 
	(alt. $q_1 - q_2 \in A$) and $q_1 - q_2 \in C$ (alt. $q_1 + q_2 \in C$). 
	
	So, there are a total of $3$ possible expressions for $\mathcal{S}^{(3)}_{q,2}$
in this case. 
	
	\item[\underline{Case 4:}] $q_1, q_2 \in A$. Then $q_1 \pm q_2 \in B$ always.
\\
	So, we get $1$ possible expression for this case. 
	\end{enumerate}
	
	Now, adding up all the possible expressions from the above four cases we get the
maximum number of possible 
	expressions for $\mathcal{S}^{(3)}_{q,2}$ to be $1 + 1 + 3 + 1 = 6$. 
\end{proof} 
With these four propositions, we are now ready for our first main theorem. 
\begin{thm} \label{theorem415}
	\begin{enumerate}
\item[(i)] Let $p \geq 5$ (alt. $p \geq 3$) be an odd prime and let $m \geq 2$ (alt.
$m \geq 3$) be any positive integer. Let 
 	$q = p^m$. Then there exist at least two $(q-6)$-dimensional orbifold lens
spaces with fundamental groups of order 
	$p^m$ which are isospectral but not isometric. 

	\item[(ii)] Let $p_1,p_2$ be odd primes such that $q = p_1 \cdot p_2 \geq 33$.
Then there exists at least two 
	$(q-6)$-dimensional orbifold lens spaces with fundamental groups of order $p_1
\cdot p_2$ which are isospectral
	but not isometric. 
	
	\item[(iii)] Let $q = 2^m$ where $m \geq 6$ be any positive integer. Then there
exist at least two $(q-5)$-dimensional 
	orbifold lens spaces with fundamental groups of order $2^m$ which are
isospectral but not isometric. 
	
	\item[(iv)] Let $q = 2p$, where $p \geq 7$ is an odd prime. Then there exist at
least two $(q-5)$-dimensional orbifold 
	lens spaces with fundamental groups of order $2p$ which are isospectral but not
isometric. 
	\end{enumerate}
\end{thm} 
\begin{proof} 
	We first recall from \propref{prop3212} that 
	\begin{align*}
	 \abs{\mathcal{L}_{0}(q,n)} \geq \sum_{t=r-2}^{r} \frac{1}{n-t}
\binom{q_{0}-1-r}{n-1-t} \binom{r}{t} 
	\end{align*}
	for $k =2$ and $r > 2$. 
	Thus for $k=2$ and $r > 2$ we have
\smallskip
	\begin{align}
	\begin{split}
	\abs{\mathcal{L}_{0}(q,n)} & \geq \frac{1}{n - (r-2)} \binom{q_0 - r -1}{n - 1
- (r-2)} \binom{r}{r-2} \\
		&\quad +\frac{1}{n-(r-1)} \binom{q_0 - r -1}{(n-1) - (r-1)} \binom{r}{r-1} +
\frac{1}{n-r} \binom{q_0 - r -1}{n-r-1}\binom{r}{r} \\
		&= \frac{1}{q_0 - 2 - r +2} \binom{q_0 - r -1}{q_0 - 2 - 1 - r +2}
\binom{r}{r-2} \\
		& \quad+ \frac{1}{q_0 -2-r+1} \binom{q_0 -r-1}{q_0 -2-1-r+1}\binom{r}{r-1} \\ 
		& \quad+ \frac{1}{q_0 - 2 -r} \binom{q_0 - r -1}{q_0 - 2 - r -1} \qquad \qquad 
\qquad  \quad \text{ since } n = q_0 - 2 \\ 
		&= \frac{1}{q_0 - r} \binom{q_0 -r -1}{q_0 - r - 1} \binom{r}{r-2} +
\frac{1}{q_0 - r -1} \binom{q_0 -r -1}{q_0 - r - 2} \binom{r}{r-1}\\
		&\quad+ \frac{1}{q_0 - r -2} \binom{q_0 - r -1}{q_0 - r -3} \binom{r}{r} \\ 
		&= \frac{1}{q_0 -r} \cdot 1 \cdot \frac{r(r-1)}{2} + \frac{1}{(q_0 - r -1)}
\cdot (q_0 -r -1) \cdot r \\
		&\quad+ \frac{1}{(q_0 -r -2)} \cdot \frac{(q_0 - r -1)(q_0 - r -2)}{2} \cdot 1
\\
		&= \frac{r(r-1)}{2(q_0 -r)} + r + \frac{(q_0 - r -1)}{2}
	\end{split}
	\end{align}
\smallskip
	It is sufficient to show that the final expression in $(3.3)$ is
greater than the number of possible expressions for the generating functions 
computed in Propositions 3.1.2--3.1.5 in order to establish the existence of 
isospectral pairs of non-isometric lens spaces. 
\smallskip
	 \begin{enumerate} 
	 \item[(i)] For $q = p^m$, we have a total of $m^2$ possible expressions for
$\tau^{(j)}_{q,2}$ from \propref{prop411}. 
	 	So, we will have isospectrality when $(3.3)$ is greater than or equal to
$m^2+1$. That is
		\begin{align}
		 & \qquad \quad \frac{r(r-1)}{2(q_0 - r)} + r + \frac{(q_0 - r - 1)}{2} \geq
m^2 +1 \notag \\
		 \Rightarrow & \quad r(r-1) + 2r(q_0 - r) +(q_0-r)(q_0-r-1) \geq 2(q_0 -
r)(m^2+1) \notag \\
		 \Rightarrow & \quad r^2 - r + (q_0 - r)[2r+q_0 - r -1 - 2m^2 - 2] \geq 0
\notag \\
		 \Rightarrow & \quad (r^2 - r) + (q_0 - r)[q_0 + r -2m^2 - 3] \geq 0 \notag 
\\
		 \Rightarrow & \quad r^2 - r + q_0^2 + q_0 r - q_0 2m^2 - 3q_0 - q_0 r - r^2 +
2m^2r + 3r \geq 0 \notag \\
		 \Rightarrow & \quad q_0(q_0 - 2m^2 - 3) + 2r(m^2 +1) \geq 0 \notag \\
		 \Rightarrow & \quad -q_0[(2m^2 + 3) - q_0] \geq -2r(m^2+1) \notag \\
		 \Rightarrow & \quad q_0[(2m^2+3) - q_0] \leq 2r(m^2+1). 
		\end{align}
		So for any given $m$, we can choose $p$ big enough so that $2m^2 + 3 \leq q_0$.
This would guarantee isospectrality. 
		We can calculate $r$ by $r = (\frac{p^{m-1}-1}{2})$ in this case. Now if $p
\geq 5$, $q_0 \geq \frac{5^m - 1}{2} > 2m^2 + 3$ for 
		all $m \geq 2$. This is easy to see since $5^m > 4m^2 + 7$ for $m \geq 2$ as the
left hand side grows exponentially greater than the right hand side. 
		So, for all $p \geq 5$ and all $m \geq 2$, $(3.4)$ will be true and we will
get isospectral pairs of dimension $(q-6) = 2n-1$. 
		 Now for $q = 3^m$, we have $3^m > 4m^2 + 7$ for $m \geq 4$. So we will have
isospectrality. We check cases $m = 2$ and $m = 3$. 
		 
		 When $m = 2$, $q = 9$, $r = 1$, $q_0 = 4$. So L.H.S. of $(3.4)$ gives
$4[2(4)+3-4] = 4(7) = 28$ and R.H.S. of $(3.4)$ gives 
		 $2(1)(4+1) = 10$. So the sufficiency condition is not satisfied. 
		 
		 When $m = 3$, $q = 27$, $r= 4$, $q_0 = 13$. L.H.S. of $(3.4)$ gives
$13[2(9)+3-13] = 13[8] = 104$ and R.H.S. of $(3.4)$ gives 
		 $2(4)[9+1] = 8(10) = 80$. So the sufficiency condition is not satisfied. 
		 
		 However, when we check individually all the possible expressions for these cases
we realize that they are less than $m^2$. 
		 
		 For $q = 3^2$, the only two expressions are for the cases when $q_1 \in A$,
$q_2 \in B_1$, $q_1 \pm q_2 \in A$, and $q_1, q_2 \in A$, $q_1 + q_2 \in A$, 
		 $q_1 - q_2 \in B_1$. No other possible expressions exist. 
		 
		 However, there are only two classes in $\mathcal{L}_0(q,n)$, i.e.,$\abs{\mathcal{L}_{0}(q,n)}= 2$. 
The two classes are
		 \begin{align*}
		 [1,2] & = \Big \{ (p_1, p_2) \in \widetilde{\mathcal{L}}_{0}(q,2) \Big
\arrowvert p_1, p_2 \in A \Big \} \\
		 [1,3] & = \Big \{ (p_1, p_2) \in \widetilde{\mathcal{L}}_{0}(q,2) \Big
\arrowvert p_1 \in A, p_2 \in B_1\text{ (alt. } p_1 \in B_1, p_2 \in A) \Big \} 
		 \end{align*} 
		 where $n = 2$, $A = \{ 1,2,4,5, 7 ,8 \}$ and $B_1 = \{ 3,6 \}$. 
		 
		 Therefore, we do not obtain isospectral pairs.  
		 
		 For $q = 3^3$, there are $7$ expressions (instead of $3^2 = 9$ possible
expressions). The case where $q_1, q_2, q_1 \pm q_2 \in B_1$ and the case where 
		 $q_1, q_2 \in B_2$ do not occur. This gives us $2$ less expressions than the
estimated number of $9$. But the number of classes is greater than or equal
		 to 
		 \begin{align*}
		  \frac{4(4-1)}{2(13-4)} + 4 + \frac{13-4-1}{2} =  \frac{2}{3} + 4 + 4 = 8
\frac{2}{3} > 7 \qquad \qquad (\text{from } (3.3)).
		 \end{align*} 
		 This means we will have non-isometric isospectral lens spaces. This gives us
our result that for $p \geq 3$ and $m \geq 3$, we will get 
		 isospectral pairs that are non-isometric. 
		 
	\item[(ii)]  For $q = p_1 \cdot p_2$, $r = \frac{p_1 + p_2 - 2}{2}$. \\
		From $(3.3)$ and \propref{prop412} we get the following sufficiency condition: 
		\begin{gather}
		\frac{r(r-1)}{2(q_0 - r)} + r + \frac{(q_0 - r -1)}{2} \geq 12 \notag \\
		\Rightarrow q_0 (25 - q_0) \leq 24 r
		\end{gather}
		From this we get that for $q_0 \geq 25$, we will always find non-isometric,
isospectral lens spaces because $(3.5)$ will always be satisfied. 
		We now check for cases where $q = 2q_0 +1 < 51$. 
		
		For $q < 51$, and $q = p_1 \cdot p_2$ with $p_1, p_2$ being odd primes, there
are only the following possibilities:
		\begin{enumerate}
			\item[(a)] $q = 3 \cdot 7 = 21$; $B = \{3,6,9,12,15,18 \}$, $C = \{7,14\}$.
\\
			In this case we have $9$ instead of $11$ possible expressions. The case where
$q_1,q_2 \in C = \{7,14\}$ is not possible, and 
			the case where $q_1, q_2 \in A$ and $q_1 \pm q_2 \in A$ is also not possible since then $q_2\equiv - q_1 (mod q)$.
Therefore, we get $2$ less expressions. now for 
			isospectrality we use $(3.3)$: 
			\begin{align*}
			\frac{4(4-1)}{2(10-4)} + 4 + \frac{(10-4-1)}{2} = 7 \frac{1}{2},
			\end{align*}
			which is not greater than $9$. So the isospectrality condition is not met. 
			
			\item[(b)] $q = 3 \cdot 5 = 15$. In this case we have $7$ instead of $11$
expressions. Here $B = \{3,6,9,12\}$ and $C = \{5,10\}$. In this 
			case, the following cases do not occur: $q_1,q_2 \in C$; $q_1 \in A, q_2 \in
C$, $q_1 \pm q_2 \in A$; $q_1, q_2, q_1 \pm q_2 \in A$; $q_1, q_2
			\in A$, $q_1 + q_2 \in A$, $q_1 - q_2 \in C$. So we get $4$ less expressions
than $11$. To check for isospectrality we use $(3.3)$ and get 
			$\frac{3(3-1)}{2(7-3)} + 3 + \frac{(7-3-1)}{2} = 5 \frac{1}{4}$, which is
less than $7$. So the isospectrality condition is not satisfied. 
			
			For (a) and (b) it can be easily shown that $\abs{\mathcal{L}_{0}(q,n)}$ is
equal to $9$ and $7$ respectively. This means that there 
			are no isospectral pairs in these cases. 
			
			\item[(c)] $q = 3 \cdot 11 = 33$. $B = \{3,6,9,12,15,18,21,24,27,30\}$ and $C
= \{11,22\}$. 
			Here $q_0 = 16$ and $r = 6$. We check for isospectrality using $(3.5)$: 
				\begin{gather*}
					q_0(25-q_0) = 16(25-16) = 144 \\
					24r = 24(6) = 144
				\end{gather*}
			So $(3.5)$ is satisfied. 
			
			\item[(d)] $q = 5 \cdot 7 = 35$, $B = \{5,10,15,20,25,30 \}$ and $C =
\{7,14,21,28 \}$.\\
			Here $q_0 = 17$ and $r = 5$. Using $(3.5)$ we get 
				\begin{gather*}
					q_0(25-q_0) = 17(25-17) = 138 \\
					24r = 24(5) = 120
				\end{gather*}
			So $(3.5)$ is not satisfied. However, we notice that in this case the actual
number of expressions is $10$ instead of $11$. So, we use 
			$(3.3)$ to check for isospectrality. Plugging in $r =5$ and $q_0 = 17$ into
$(3.3)$ we get 
				\begin{align*}
					\frac{5(4)}{2(12)} + 5 + \frac{11}{2} = 11 \frac{1}{3} > 10
				\end{align*}
			This implies that $\mathcal{S}^{(3)}_{q,2}$ is not one-one and therefore, we
will have non-isometric isospectral lens spaces in this case. 
			
			\item[(e)] Finally, we check $q = 3 \cdot 13 = 39$. \\
			Here $q_0 = 19$ and $r = 7$. Using $(3.5)$ we see 
				\begin{gather*}
					q_0(25-q_0) = 19(25-19) = 114 \\
					24r = 24(7) = 168
				\end{gather*}
			So $(3.5)$ is satisfied and we will have isospectral pairs in this case. \\
			(a)-(e) are all the cases of $q = p_1 \cdot p_2 < 51$, where $p_1, p_2$ are
odd primes. \\ \\
			Combining these results with the fact that for $q \geq 51$, $(3.5)$ will
always be satisfied, we prove (iii). 
		\end{enumerate}
		
		\item[(iii)] Let $q = 2^m$. We use \propref{prop413} and $(3.3)$ to get a
sufficiency condition for isospectrality:
			\begin{align*}
				\frac{r(r-1)}{2(q_0 - r)} + r + \frac{(q_0 - r -1)}{2} \geq (m-1)^2 + 1
			\end{align*}
		Here $q_0 = \frac{2^m}{2} = 2^{m-1}$ and $2r = 2^{m-1}$. Therefore, $q_0 = 2r$
in this case. Simplifying the above inequality, 
		we get 
			\begin{align*}
				q_0 [(2m^2 - 4m +5) - q_0] \leq 2r (m^2 - 2m + 2).
			\end{align*}
		But since $q_0 = 2r$, we get 
			\begin{align}
				(m^2 - 2m + 3) \leq q_0
			\end{align}
		If $m \geq 6$, then $m^2 - 2m + 3 < 2^{m-1}$. Further, the right hand side of
$(3.6)$ grows exponentially bigger than the left hand side 
		as $m$ grows. For $m = 3$, $4$ and $5$, the actual number of expressions for
$\tau_{q,2}^{(j)}$ are $4$, $9$ and $16$ respectively. 
		Further, it can be easily shown that for $m= 3$, $4$ and $5$, $\abs{
\mathcal{L}_{0}(q,n)}$ is $4$, $9$ and $16$ respectively. 
		Therefore, for $m=3$, $4$ and $5$ we do not get isospectrality. This gives us
(iii). 
		
		\item[(iv)] Using \propref{prop414} and $(3.3)$ we get the sufficiency
condition for isospectrality for $q = 2p$, where $p$ is an odd prime 
		$\geq 7$. Note that in this case $q_0 = \frac{q}{2} = p$ and $r =
\frac{q+2}{4}$. Now for isospectrality we should have 
			\begin{align*}
				& \qquad \frac{r(r-1)}{2(q_0 - r)} + r + \frac{(q_0 - r -1)}{2} \geq 7 \\
				\Rightarrow \, \, & q_0(15 - q_0) \leq 14r \\
				\Rightarrow \, \, & p(15-p) \leq 7(p+1) \\
				\Rightarrow \, \, & 0 \leq p^2 - 8p + 7 \\
				\text{or } & (p-1)(p-7) \geq 0
			\end{align*}
		Since $p$ is positive, whenever $p \geq 7$, we will have isospectrality. When
$ q = 2 \cdot 5= 10$, then $\abs{ \mathcal{L}_{0}(q,n)} = 6 = 
		$ number of expressions for $\mathcal{S}_{q,2}^{(3)}$. So, we do not get
isospectral pairs when $p=5$. This proves (iv). 
	\end{enumerate}
\end{proof}

\subsection{Lens Spaces for General Integers} \label{subsection43}

\begin{align*}
	\text{Let }\quad \,  L &= L(q: p_1, \ldots, p_n) = S^{2n-1} / G \qquad
\text{and} \\
			L' &= L(q: p_1, \ldots, p_n) = S^{2n-1} / G' 
\end{align*}
be two isospectral non-isometric orbifold lens spaces as obtained in
\secref{section41} where $G = \langle g \rangle$, $G' = \langle g' \rangle$. 
\[ g = 
\begin{pmatrix}
 R(p_{1} / q) & & \text{ {\huge 0}} \\
  & \ddots & \\
  \text{ {\huge 0}} & & R(p_{n} / q)
\end{pmatrix}
\]
and 
\[ g' = 
\begin{pmatrix}
 R(s_{1} / q) & & \text{ {\huge 0}} \\
  & \ddots & \\
  \text{ {\huge 0}} & & R(s_{n} / q)
\end{pmatrix}
\]
We define 
\[ \tilde{g}_{W+} = 
\begin{pmatrix}
 R(p_{1} / q) & & & \text{ {\huge 0}} \\
  & \ddots & & \\
  & & R(p_{n} / q)& \\
   \text{ {\huge 0}} & & & I_W
\end{pmatrix}
\]
and
\[ \tilde{g}'_{W+} = 
\begin{pmatrix}
 R(s_{1} / q) & & & \text{ {\huge 0}} \\
  & \ddots & & \\
  & & R(s_{n} / q)& \\
   \text{ {\huge 0}} & & & I_W
\end{pmatrix}
\] \\
where $I_W$ is the $W \times W$ identity matrix for some integer $W$. We can
define $\tilde{G}_{W+}$ $= \langle \tilde{g}_{W+} \rangle$ and 
$\tilde{G}'_{W+} = \langle \tilde{g}'_{W+} \rangle$. Then $\tilde{G}_{W+}$ and
$\tilde{G}'_{W+}$ are cyclic groups of order $q$. We define lens spaces 
$\tilde{L}_{W+} = S^{2n+W-1} / \tilde{G}_{W+}$ and $\tilde{L}'_{W+} = S^{2n+W-1}
/ \tilde{G}'_{W+}$. Then, like \propref{lensspaceisometric}, we get: 

\begin{prop}\label{prop431}
	Let $\tilde{L}_{W+}$ and $\tilde{L}'_{W+}$ be as defined above. Then
$\tilde{L}_{W+}$ is isometric to $\tilde{L}'_{W+}$ iff there is a number $l$ coprime 
with $q$ and there are numbers $e_{i} \in \{-1, 1\}$ such that $(p_1, \ldots, p_n)$ is a
permutation of $(e_1 l s_1, \ldots, e_n l s_n)\pmod{q}$. 
\end{prop}

The following lemma follows directly from this proposition. 

\begin{lemma}\label{lemma432}	
	Let $L$, $L'$, $\tilde{L}_{W+}$ and $\tilde{L}'_{W+}$ be as defined above. Then
$L$ is isometric to $L'$ iff $\tilde{L}_{W+}$ is isometric to 
	$\tilde{L}'_{W+}$.
\end{lemma}

We get the following theorem: 

\begin{thm}\label{theorem433}
	Let $\widetilde{\mathcal{L}}_{0}^{W+}(q,n,0)$ be the family of all $(2n+W-1)$-dimensional
orbifold lens spaces with fundamental groups of order $q$
	that are obtained in the manner described above. Let $\tilde{L}_{W+} \in
\mathcal{L}_{0}^{W+}(q,n,0)$ (where $\mathcal{L}_{0}^{W+}(q,n,0)$
	denotes the set of isometry classes of
$\widetilde{\mathcal{L}}_{0}^{W+}(q,n,0)$). Let $F_{q}^{W+}(z: p_1, \ldots, p_n,
0)$ be the generating 
	function associated to the spectrum of $\tilde{L}_{W+}$. Then on the domain
$\big \{z \in \mathbf{C} \big \arrowvert \abs{z} < 1 \big \}$, 
	\begin{align*}
		F_{q}^{W+}(z: p_1, \ldots, p_n, 0) = \frac{ (1 + z )}{(1-z)^{W-1}} \cdot
\frac{1}{q} \sum_{l=1}^{q} \frac{1}{\prod_{i = 1}^{n} (z - \gamma^{p_i l})(z -
\gamma^{-p_i l})}
	\end{align*}
\end{thm}
\begin{proof}

Recall the definitions of $\Delta_0$, $r^2$, $P^k$, $H^k$, $\mathcal{H}^k$ and
$\mathcal{H}^{k}_{G}$ from \secref{section32}. We extend the definitions for 
	$\mathbf{R}^{2n+W}$. That is, let $\Delta_0$ be the Laplacian on
$\mathbf{R}^{2n+W}$ with respect to the Flat Riemannian metric;
	$r^2 = \sum_{i=1}^{2n+W} x_{i}^2$, where $(x_1, \ldots, x_{2n+W})$ is the
standard coordinate system on $\mathbf{R}^{2n+W}$; 
	for $k \geq 0$, $P^k$ is the space of complex valued homogenous polynomials of
degree $k$ in $\mathbf{R}^{2n+W}$; $H^k$ is the subspace of 
	$P^k$ consisting of harmonic polynomials on $\mathbf{R}^{2n+W}$;
$\mathcal{H}^k$ is the image of $i^{*}: C^{\infty}(\mathbf{R}^{2n+W})
\longrightarrow C^{\infty}
	(S^{2n+W-1})$ where $i: S^{2n+W-1} \longrightarrow \mathbf{R}^{2n+W}$ is the natural
injection; and $\mathcal{H}_{\tilde{G}}^{k}$ is the space of all
$\widetilde{G}$-invariant 
	functions of $\mathcal{H}^k$.
	
	Then from \propref{ch3inducedlinmap} and \propref{ch3eigenspaceprop}, we get
that $H^k$ is $O(2n+W)$-invariant; $P^k$ has the direct sum decomposition 
	$P^k = H^k \oplus r^2 P^{k-2}$; $\mathcal{H}^k$ is an eigenspace of
$\widetilde{\Delta}$ on $S^{2n+W-1}$ with eigenvalue $k(k + 2n + W -2)$; 
	$\sum_{k=0}^{\infty} \mathcal{H}^k$ is dense in $C^{\infty}(S^{2n+W-1})$ in the
uniform convergence topology and $\mathcal{H}^{k}$ is isomorphic to $H^k$. 
	
	This along with the results in \corref{corollary327}, where $\dim E_{k(k+2n + W-1)}
= \dim \mathcal{H}_{\tilde{G}_{W+}}^{k}$, we get 
	\begin{align*}
		 F_{q}(z: p_1, \ldots, p_n, 0) = \sum_{k=0}^{\infty} (\dim
\mathcal{H}_{\tilde{G}_{W+}}^{k}) z^k.
	\end{align*}
	Now $\tilde{G}_{W+}$ is contained in $SO(2n+W)$. 
	
	Let $\chi_k$ and $\tilde{\chi}_k$ be the characters of the natural
representations of $SO(2n+W)$ on $H^k$ and $P^k$ respectively. Then 
	\begin{align}
	\dim \mathcal{H}_{\widetilde{G}_{W+}}^{k} = \frac{1}{\abs{\widetilde{G}_{W+}}}
\sum_{\tilde{g}_{W+} \in \widetilde{G}} \chi_k (\tilde{g}_{W+}) = \frac{1}{q}
\sum_{l = 1}^{q} \chi_k (\tilde{g}_{W+}^l)
	\end{align}
	\propref{ch3inducedlinmap} gives 
	\begin{align}
		\chi_{k}(\tilde{g}_{W+}^l) = \tilde{\chi}_{k}(\tilde{g}_{W+}^l) -
\tilde{\chi}_{k-2}(\tilde{g}_{W+}^l)
	\end{align}
	where $\tilde{\chi}_t = 0$ for $t > 0$. 
	
If $W$ is even, then we can view the space $P^k$ as having a basis consisting of all monomials of
the form: 
\begin{align*}
	z^{I} \cdot \bar{z}^{J} = (z_1)^{i_1} \cdots (z_{n+v})^{i_{n+v}} \cdot
(\bar{z}_{1})^{j_1} \cdots (\bar{z}_{n+v})^{j_{n+v}} 
\end{align*}
where $W = 2v$ and where 
$I_{n+v} + J_{n+v} = i_1 + \cdots + i_{n+v} + j_1 + \cdots + j_{n+v} = k$
and $i_1 , \ldots, i_{n+v}, j_1, \ldots, j_{n+v} \geq 0$. Then, 
\begin{align*}
	\tilde{g}_{W+}^{l} (z^{I} \cdot \bar{z}^{J}) = \gamma^{i_1 p_1 l + \cdots + i_n
p_n l - j_1 p_1 l - \cdots - j_n p_n l} (z^{I} \cdot \bar{z}^{J}).
\end{align*}
If $W$ is odd, (say $W = 2u +1$), then we get for basis of $P^k$ 
\begin{align*}
	z^{I} \cdot \bar{z}^{J} \cdot z^{t}_{n + 2u + 1} = (z_1)^{i_1} \cdots
(z_{n+u})^{i_{n+u}} \cdot (\bar{z}_1 )^{j_1} \cdots (\bar{z}_{n+u})^{j_{n+u}}
\cdot (z_{n + 2u + 1})^{t}
\end{align*}
where $z_{n + 2u + 1} = x_{n+W}$ where $(x_1, y_1, \ldots, x_{n+W - 1}, y_{n+W -
1}, x_{n + W})$ is the standard euclidean coordinate system on
$\mathbf{R}^{2n+W}$ 
with $z_i = x_i + iy_i$ for $i = 1, 2, \ldots, n+W - 1$, and $i_1, \ldots,
i_{n+u}, j_1, \ldots, j_{n+u}$, $t \geq 0$ and $i_1 + \cdots + i_{n+u} + j_1 +
\cdots + j_{n+u} + t = k 
= I_{n + u} + J_{n+u} + t$. So, in that case
\begin{align*}
	\tilde{g}_{W+}^{l} (z^{I} \cdot \bar{z}^{J} \cdot z_{n + 2u +1} ) = \gamma^{i_1
p_1 l + \cdots + i_n p_n l - j_1 p_1 l - \cdots - j_n p_n l} (z^{I} \cdot
\bar{z}^{J} \cdot z_{n + 2u +1} ) 
\end{align*}
So, for $W$, even case, we will get 
\begin{gather*}
	\begin{split}
	F_{q}^{W+}(z: p_1, \ldots, p_n, 0) &= \frac{1}{q} \sum_{k=0}^{\infty} \, \sum
_{l =1}^{q} \chi_{k} ( \tilde{g}_{W+}^{l} ) z^k \\
		&=\frac{ (1-z^2) }{q} \sum_{l =1}^{q} \, \sum _{k=0}^{\infty} \tilde{\chi}_{k}
(\tilde{g}_{W+}^{l} )z^k \\
		&=\frac{ (1-z^2) }{q} \sum_{l=1}^{q} \,  \sum_{k=0}^{\infty} \sum_{I_{n+v} +
J_{n+v} = k} \gamma^{i_1 p_1 l + \cdots + i_np_nl - j_1 p_1 l - \cdots - j_n p_n
l}z^k
	\end{split}
\end{gather*}
\begin{multline*}
	=\frac{ (1-z^2) }{q} \sum_{l=1}^{q} \, \sum_{k=0}^{\infty} \sum_{I_{n+v} +
J_{n+v} = k} (\gamma^{p_1 l} z)^{i_1} \cdots  (\gamma^{p_n l}
z)^{i_n}(\gamma^{-p_1 l} z)^{j_1} \cdots  \\
		 (\gamma^{-p_n l} z)^{j_n} \cdot z^{i_{n+1} + \cdots +  i_{n+v} + j_{n+1} +
\cdots + j_{n+v}}
\end{multline*}
\begin{align*}
	&=\frac{(1-z^2)}{q} \sum_{l=1}^{q} \, \prod_{i=1}^{n} (1 + \gamma^{p_i l}z +
\gamma^{2p_i l}z^2 + \cdots)(1 + \gamma^{-p_i l}z + \gamma^{-2p_i l}z^2 +
\cdots) 
		 (1+z + z^2 + \cdots)^W \\
	&=\frac{(1-z^2)}{q} \sum_{l=1}^{q} \frac{1}{\prod_{i=1}^{n} (1 - \gamma^{p_i
l}z)(1- \gamma^{-p_i l}z)(1-z)^W} \quad \text{on } \big \{z \in \mathbf{C} \big
\arrowvert \abs{z} < 1\big \} \\ 
	&= \frac{(1+z)}{(1-z)^{W-1}} \cdot \frac{1}{q} \sum_{l=1}^{q}
\frac{1}{\prod_{i=1}^{n} (z - \gamma^{p_i l})(z- \gamma^{-p_i l})}
\end{align*}
For $W$ odd case, we get by similar calculations, 
\begin{align*}
	F_{q}^{W+}(z: p_1, \ldots, p_n) &= \frac{ (1-z^2) }{q} \sum_{l=1}^{q} \, 
\sum_{k=0}^{\infty} \sum_{I_{n+u} + J_{n+u} + t = k} \gamma^{i_1 p_1 l + \cdots
+ i_np_nl - j_1 p_1 l - \cdots - j_n p_n l}z^k
\end{align*}
\begin{multline*}
	=\frac{ (1-z^2) }{q} \sum_{l=1}^{q} \, \sum_{k=0}^{\infty} \sum_{I_{n+u} +
J_{n+u}+t = k} (\gamma^{p_1 l} z)^{i_1} \cdots  (\gamma^{p_n l}
z)^{i_n}(\gamma^{-p_1 l} z)^{j_1} \cdots  \\
		 (\gamma^{-p_n l} z)^{j_n} \cdot z^{i_{n+1} + \cdots +  i_{n+u} + j_{n+1} +
\cdots + j_{n+u}+ t}
\end{multline*}
\begin{align*}
	&=\frac{(1-z^2)}{q} \sum_{l=1}^{q} \, \prod_{i=1}^{n} (1 + \gamma^{p_i l}z +
\gamma^{2p_i l}z^2 + \cdots)(1 + \gamma^{-p_i l}z + \gamma^{-2p_i l}z^2 +
\cdots) 
		 (1+z + z^2 + \cdots)^W \\
	&=\frac{(1-z^2)}{q} \sum_{l=1}^{q} \frac{1}{\prod_{i=1}^{n} (1 - \gamma^{p_i
l}z)(1- \gamma^{-p_i l}z)(1-z)^W} \quad \text{on } \big \{z \in \mathbf{C} \big
\arrowvert \abs{z} < 1\big \} \\ 
	&= \frac{(1+z)}{(1-z)^{W-1}} \cdot \frac{1}{q} \sum_{l=1}^{q}
\frac{1}{\prod_{i=1}^{n} (z - \gamma^{p_i l})(z- \gamma^{-p_i l})} \quad \text{
as before.} 
\end{align*}
\end{proof}

\begin{cor}\label{cor434} When $L(q: p_1, \ldots, p_n)$ and $L(q: s_1, \ldots, s_n)$ have the same generating function, then $\tilde{L}_{W+}$ and $\tilde{L}'_{W+}$ (as defined above) also have the same generating function
\end{cor}

\begin{proof} This follows from the fact that \[{{F_{q}}^{W+}}(z: p_1, \ldots, p_n, 0) = \frac{1}{{(1-z)}^W}\ F_{q}(z: p_1, \ldots, p_n)\].
\end{proof} 

The above results give us the following theorem. 
\begin{thm}\label{theorem435} 
\begin{enumerate}
\item[(i)] Let $P \geq 5$ (alt. $P \geq 3$) be any odd prime and let $m \geq 2$ (alt.
$m \geq 3$) 
be any positive integer. Let $q =P^{m}$. Then there exist at least two
$(q+W-6)$-dimensional orbifold lens spaces with fundamental groups of 
order $P^{m}$ which are isospectral but not isometric. 
\item[(ii)] Let $P_{1},P_{2}$ be two odd primes such that $q=P_{1}\cdot P_{2}
\geq 33$. Then there exist at least two $(q+W-6)$-dimensional orbifold 
lens spaces with fundamental groups of order $P_{1} \cdot P_{2}$ which are
isospectral but not isometric. 
\item[(iii)] Let $q = 2^{m}$ where $m \geq 6$ is any positive integer. Then
there exist at least two $(q+W-5)$-dimensional orbifold lens spaces with
fundamental groups of order $2^{m}$ which are isospectral but not isometric. 
\item[(iv)] Let $q=2P$, where $P \geq 7$ is an odd prime. Then there exist at
least two $(q+W-5)$-dimensional orbifold lens spaces with 
fundamental groups of order $2P$ which are isospectral but not isometric. 
\end{enumerate}
\end{thm}

\begin{cor}\label{corollary436}
\begin{enumerate}
\item[(i)] Let $x \geq 19$ be any integer. Then there exist at least two
$x$-dimensional orbifold lens spaces with fundamental groups of order $25$ which
are isospectral but not isometric. 
\item[(ii)] Let $x \geq 27$ be any integer. Then there exist at least two
$x$-dimensional orbifold lens spaces with fundamental group of order $33$ which
are isospectral but not isometric. 
\item[(iii)] Let $x \geq 59$ be any integer. Then there exist at least two $x$
dimensional orbifold lens spaces with fundamental group of order $64$ which are
isospectral but not isometric. 
\item[(iv)] Let $x \geq 9$ be any integer. Then there exist at least two $x$
dimensional orbifold lens spaces with fundamental group of order $14$ which are
isospectral but not isometric. 
\end{enumerate}
\end{cor}

\begin{proof}
(i) Let $q = 25$ and $W \in \{0, 1, 2, 3, \ldots \}$ in $(i)$ of the theorem. 
\begin{enumerate}
	\item[(ii)] Let $q =33$ and $W \in \{0, 1, 2, 3, \ldots \}$ in $(ii)$ of the
theorem. 
	\item[(iii)] Let $q =64$ and $W \in \{0, 1, 2, 3, \ldots \}$ in $(iii)$ of the
theorem. 
	\item[(iv)] Let $q =14$ and $W \in \{0, 1, 2, 3, \ldots \}$ in $(iv)$ of the
theorem. 
\end{enumerate}
\end{proof}

When W is an odd numnber, we get even dimensional orbifold lens spaces that are 
isospectral but not isometric.





\section{Examples}\label{examples}
In this section we will look at some examples of isospectral
non-isometric orbifold lens spaces by calculating their 
generating functions. We will also look at an example that will suggest that our
technique can be generalized for higher values
of $k = q_0 - n$. Recall that $q_0 = (q-1)/2$ for odd $q$ and $q_0 = q / 2$ for
even $q$. In the previous sections we assumed 
$k = 2$. The technique for getting examples for higher values of $k$ is similar,
but as we shall see, the calculations for the different
types of generating functions becomes more difficult as $k$ increases. 

In all the examples we will denote a lens space by $L(q: p_1, \ldots, p_n) =
S^{2n-1} / G$, where $G$ is the cyclic group
generated by $g =
\begin{pmatrix}
 R(p_{1} / q) & & \text{ {\huge 0}} \\
  & \ddots & \\
  \text{ {\huge 0}} & & R(p_{n} / q)
\end{pmatrix} $. We will write $G = \langle g \rangle$. 

\subsection{Examples for $k =2$} \label{section51}

\begin{ex} \label{example512}
Let $q = 5^2 = 25$, $q_0 = \frac{q-1}{2} = 12$, $n = 10$, $k =2$, \\ $A = \{ 1,
2, 3, 4 , 6 ,7,8,9,11, 12 ,13, 14, 16 , 17, 18, 19, 21 , 22, 23, 24 \}$,
$B_1 = \{ 5, 10, 15, 20\}$. Let $w([p_1, \ldots, p_{10}]) = [q_1, q_2]$. Let
$\gamma = e^{2 \pi i / 25}$ and $\lambda = e^{2 \pi i / 5}$. $a_0 = \abs{A} =
20$, $b_{0,1} = \abs{B_1} = 4$, 
$\sum_{l \in A} \gamma^l = 0$, $\sum_{l \in B_1} \gamma^l = -1$, $\sum_{l \in A}
\lambda^l = -5$ and $\sum_{l \in B_1} \lambda^l = 4$ (from (4.1)). 
\begin{enumerate}
\item[\underline{Case 1:}] $q_1, q_2 \in B_1$ and $q_1 \pm q_2 \in B_1$. So,  
	\begin{align*}
		\psi_{25,2}([q_1, q_2])(z) &= 20z^4 + 20z^3 + 20z^2 + 20z + 20 \\
		\alpha_{25,2}^{(1)} ([q_1 , q_2])(z) & = 4z^4 - 16z^3 + 24z^2 - 16z + 4
	\end{align*}
This corresponds to the case where $[p_1, \ldots, p_{10}] =
[1,2,3,4,6,7,8,9,11,12]$ 
which corresponds to a manifold lens spaces. 

\item[\underline{Case 2:}] Since there is only one $B_1$ this case does not
occur. 
	
\item[\underline{Case 3:}] $q_1 \in B_1$ and $q_2 \in A$ $(\text{alt. } q_1 \in
A, q_2 \in B_1)$. $q_1 \pm q_2 \in A$ always. \\
	So,  \begin{align*}
		\psi_{25,2}([q_1, q_2])(z) &= 20z^4 + 10z^3 + 40z^2 + 10z + 20 \\
		\alpha_{25,2}^{(1)} ([q_1 , q_2])(z) &= 4z^4 - 6z^3 + 4z^2 - 6z + 4 \\
		\text{corresponding to} \quad [p_1, \ldots, p_{10}]& = [1,2,3,4,5,6,7,8,9,11]
\\
		\text{and to} \quad  [s_1, \ldots, s_{10}] &= [1,2,3,4,6,7,8,9,10,11] \\ 
		\text{and} \quad [p_1, \ldots, p_{10}]  & \neq [s_1, \ldots, s_{10}]. 
	\end{align*}

So, we get two isospectral non-isometric orbifolds: $L_{1} = L(25: 1, 2, 3, 4,
5, 6, 7, 8 ,9 ,11)$ and $L_2 = (25: 1, 2, 3, 4, 6, 7, 8, 9, 10, 11)$. 
We denote by $\sum_{i}$ the singular set of $L_i$. Then
	$\sum_{1} = \big \{ (0,0,\ldots , x_9, x_{10}, 0, 0, \ldots, 0) \in S^{19} \big
\arrowvert x_{9}^2 + x_{10}^{2} = 1 \big \}$
	and  $ \sum_{2} = \big \{ (0,0, \ldots, x_{17}, x_{18}, 0, 0) \in S^{19} \big
\arrowvert x_{17}^{2} + x_{18}^{2} = 1 \big \}$ with 
	isotropy groups $\langle g_{1}^5 \rangle$ and $\langle g_{2}^{5} \rangle $
where
\[ g_{1}^5 = 
\begin{pmatrix}
 R(5p_{1} / 25) & & \text{ {\huge 0}} \\
  & \ddots & \\
  \text{ {\huge 0}} & & R(5p_{10} / 25)
\end{pmatrix}
= 
\begin{pmatrix}
 R(p_{1} / 5) & & \text{ {\huge 0}} \\
  & \ddots & \\
  \text{ {\huge 0}} & & R(p_{10} / 5)
\end{pmatrix}
\]
and 
\[ g_{2}^5 = 
\begin{pmatrix}
 R(s_{1} / 5) & & \text{ {\huge 0}} \\
  & \ddots & \\
  \text{ {\huge 0}} & & R(s_{10} / 5)
\end{pmatrix}
\]
where $g_1$ and $g_2$ are generators of $G_1$ and $G_2$, respectively with $L_1
= S^{19} / G_1$ and $L_2 = S^{19} / G_2$. 
$\sum_{1}$ and $\sum_{2}$ are homeomorphic to $S^1$. We denote the two isotropy
groups by $H_1 = \langle g_{1}^5 \rangle$ and 
$H_2 = \langle g_{2}^{5} \rangle$. 

\item[\underline{Case 4:}] $(a)$ $q_1, q_2 \in A$ and $q_1 \pm q_2 \in A$. So, 
	\begin{align*}
		\psi_{25,2}([q_1, q_2])(z) &= 20z^4 + 40z^2 + 20 \\
		\alpha_{25,2}^{(1)} ([q_1 , q_2])(z) &= 4z^4 + 4z^3 + 4z^2 + 4z + 4
	\end{align*}
corresponding to 
	\begin{align*}
		L_3 &= L(25: 1, 2, 3, 4, 5, 6 ,7, 8, 9, 10) = S^{19} / G_{3}, \text{ where }
G_{3} = \langle g_3 \rangle \\ 
		L_4 &= L(25: 1, 2, 3, 4, 5, 6 ,7, 8, 10, 11) = S^{19} / G_{4}, \text{ where }
G_{4} = \langle g_4 \rangle \\ 
	\text{and} \quad L_5 &= L(25: 1, 2, 3, 4, 5, 6 ,7, 10, 11, 12) = S^{19} /
G_{5}, \text{ where } G_{5} = \langle g_5 \rangle
	\end{align*}
The isotropy groups for $L_3$, $L_4$ and $L_5$ are $\langle g_{3}^5 \rangle$,
$\langle g_{4}^{5} \rangle$ and $\langle g_{5}^{5} \rangle$, 
respectively. $\sum_{3}$, $\sum_{4}$ and $\sum_{5}$ are all homeomorphic to
$S^3$. So, here we get $3$ isospectral orbifold lens spaces 
that are non-isometric. \\

$(b)$ $q_1, q_2 \in A$ and $q_1 + q_2 \in B_1$, $q_1 - q_2 \in A$ $(\text{alt. }
q_1 + q_2 \in A, \, q_1 - q_2 \in B_1)$. So, 
	\begin{align*}
		\psi_{25,2}([q_1, q_2])(z) &= 20z^4 + 30z^2 + 20 \\
		\alpha_{25,2}^{(1)} ([q_1 , q_2])(z) &= 4z^4 + 4z^3 + 14z^2 + 4z + 4
	\end{align*}
corresponding to 
	\begin{align*}
		L_6 &= L(25: 1, 2, 3, 4, 5, 6 ,7, 9, 10, 11) = S^{19} / G_{6}, \text{ where }
G_{6} = \langle g_6 \rangle \\ 
		\text{and} \quad L_7 &= L(25: 1, 2, 3, 4, 5, 6 ,7, 8, 10, 11) = S^{19} /
G_{7}, \text{ where } G_{7} = \langle g_7 \rangle
	\end{align*}
Then, again, $\sum_{6}$ and $\sum_{7}$ are homeomorphic to $S^3$, and $L_6$ and
$L_7$ have isotropy groups $\langle g_{6}^5 \rangle$ 
and $\langle g_{7}^5 \rangle$. 
\end{enumerate}
\end{ex}

\begin{ex} \label{example513}
$q = 3^3 = 27$; $q_0 = 13$, $k = 2$, $n = 11$ and \\ 
$A = \{ 1,2, 4, 5, 7, 8, 10, 11 ,13 , 14, 16, 17, 19, 20, 22, 23, 25, 26\}$,
$B_1 = \{3, 6, 12, 15, 21, 24 \}$, $B_2 = \{9, 18 \}$. 
Let $w([p_1, \ldots, p_{11}]) = [q_1, q_2]$. Let $\gamma = e^{2 \pi i /27}$,
$\lambda = e^{2 \pi i / 9}$ and $\delta = e^{2 \pi i / 3}$ 
be primitive $27^{\text{th}}$, $9^{\text{th}}$ and $3^{\text{rd}}$ roots of
unity, respectively. Here we get isospectral non-isometric pairs only in two cases:
\begin{enumerate}
\item[\underline{Case 1:}] $q_1 \in B_1$ and $q_2 \in A$ $(\text{alt. } q_1 \in A \text{ and } q_2
\in B_1 )$. $q_1 \pm q_2 \in A$ always. So we get, 
		\begin{align*} 
			\psi_{27,2}([q_1, q_2])(z) &= 18z^4 + 36z^2 + 18 \\
			\alpha_{27,2}^{(1)} ([q_1 , q_2])(z) &= 6z^4 + 6z^3 +12z^2 + 6z + 6 \\ 
			\alpha_{27,2}^{(2)} ([q_1 , q_2])(z) &= 2z^4 - 2z^3 - 2z + 2		
		\end{align*}
	corresponding to 
		\begin{align*}
			L_1 &= L(27: 1, 2, 4, 5, 6, 7, 9, 10, 11, 12 ,13), \\ 
			L_2 &= L(27: 1, 4, 5, 6, 7, 8, 9, 10, 11, 12 ,13) \quad \text{ and } \\ 
			L_3 &= L(27: 1, 2, 5, 6, 7, 8, 9, 10, 11, 12, 13) 
		\end{align*}
	If $G_1 = \langle g_1 \rangle$, $G_2 = \langle g_2 \rangle$, $G_3 = \langle g_3
\rangle$ are such that  $L_1 = S^{21} / G_1$, $L_2 = S^{21} / G_{2}$ and 
	$L_{3} = S^{21} / G_{3}$, then 
		\begin{align*}
		\textstyle \sum_{1} &= \big \{ (0,0, \ldots, x_9 , x_{10}, 0, x_{13}, x_{14} ,
0 , \ldots,  x_{19} , x_{20} , 0, 0 ) \in S^{21} \text{: isotropy group}
			 = \langle g_{1}^{9} \rangle \big \} \\ 
		& \bigcup \big \{ (0,0, \ldots, x_{13}, x_{14}, 0 , \ldots, 0) \in S^{21}
\text{: isotropy group} = \langle g_{1}^{3} \rangle \big \} \\ 
		 \textstyle  \sum_{2} &= \big \{ (0,0, \ldots, x_{7}, x_{8}, 0, \ldots,
x_{13}, x_{14}, \ldots, 0, x_{19}, x_{20}, 0 ,0 ) \in S^{21} \text{: isotropy
group}
		 	 = \langle g_2^9 \rangle \big \} \\ 
		& \bigcup \big \{ (0,0, \ldots, x_{13}, x_{14}, 0, \ldots, 0) \in S^{21}
\text{: isotropy group} = \langle g_2^3 \rangle \big \} \\
		\textstyle \sum_{2} &= \big \{ (0,0, \ldots, x_7, x_8, 0, \ldots, x_{13},
x_{14}, \ldots, x_{19}, x_{20}, 0 ,0) \in S^{21} \text{: isotropy group}
			= \langle g_3^9 \rangle \big \} \\ 
		& \bigcup \big \{ (0,0, \ldots, x_{13}, x_{14}, 0 , \ldots, 0) \in S^{21}
\text{: isotropy group} = \langle g_3^3 \rangle \big \} 
		\end{align*}
	So all three orbifolds have the same isotropy type and all the singular sets are homeomorphic to
$S^5$. \\
	
 \item[\underline{Case 2:}] $q_1 + q_2 \in B_1$, $q_1 - q_2 \in A$ $(\text{alt. } q_1 + q_2 \in A \,
, \, q_1 - q_2 \in B_1 )$. So we get, 
		\begin{align*} 
			\psi_{27,2}([q_1, q_2])(z) &= 18z^4  +  36z^2  + 18 \\
			\alpha_{27,2}^{(1)} ([q_1 , q_2])(z) &= 6z^4 + 6z^2 + 6 \\ 
			\alpha_{27,2}^{(2)} ([q_1 , q_2])(z) &= 2z^4 + 4z^3 + 6z^2 + 4z + 2		
		\end{align*}
	corresponding to
		\begin{align*}
			L_{4} &= L(27: 1, 3, 4, 5, 6, 7, 9, 10, 11, 12, 13) = S^{21} / G_4 \, ; \,
G_4 = \langle g_4 \rangle \\ 
			L_{5} &= L(27: 1, 3, 4, 5, 6, 7, 8, 9, 11, 12, 13) = S^{21} / G_5 \, ; \, G_5
= \langle g_5 \rangle \\ 
			L_{6} &= L(27: 1, 3, 5, 6, 7, 8, 9, 10, 11, 12, 13) = S^{21} / G_6 \, ; \,
G_6 = \langle g_6 \rangle
		\end{align*}
	Then, 
	\begin{gather*}
	\begin{split}
		\textstyle \sum_{4} &= \big \{ (0,0,x_3, x_4, 0,\ldots, 0, x_9, x_{10}, 0,0,
x_{13}, x_{14}, 0, \ldots, x_{19}, x_{20}, 0, 0) \in S^{21} \\
			 &\qquad \qquad \qquad \qquad \qquad \qquad \qquad \qquad \qquad \text{ with
isotropy group} = \langle g_{4}^{9} \rangle \big \} \\ 
		& \bigcup \big \{ (0,0, \ldots,0,  x_{13}, x_{14}, 0 , \ldots, 0) \in S^{21} 
\text{ with isotropy group} = \langle g_{4}^{3} \rangle \big \} \\ \\
		 \textstyle  \sum_{5} &= \big \{ (0,0,x_3,x_4,0,\ldots, 0, x_9, x_{10}, 0,
\ldots, 0, x_{15}, x_{16}, 0, 0, x_{19}, x_{20}, 0, 0) \in S^{21}\\
		  	&\qquad \qquad \qquad \qquad \qquad \qquad \qquad \qquad \qquad \text{ with
isotropy group} = \langle g_5^9 \rangle \big \} \\ 
		& \bigcup \big \{ (0,0, \ldots, 0, x_{15}, x_{16}, 0, \ldots, 0) \in S^{21}
\text{ with isotropy group} = \langle g_5^3 \rangle \big \}
	\end{split}
	\end{gather*}
	\begin{align*}
		\textstyle \sum_{6} &= \big \{ (0,0,x_3, x_4, 0, 0, x_7, x_8, 0, \ldots, 0,
x_{13}, x_{14}, 0, \ldots, 0, x_{19}, x_{20}, 0, 0) \in S^{21} \\ 
			&\qquad \qquad \qquad \qquad \qquad \qquad \qquad \qquad \qquad \text{ with
isotropy group} = \langle g_6^9 \rangle \big \} \\ 
		& \bigcup \big \{ (0,0, \ldots, 0,  x_{13}, x_{14}, 0 , \ldots, 0) \in S^{21}
\text{ with isotropy group} = \langle g_6^3 \rangle \big \} 
	\end{align*}
	So, the singular sets are all homeomorphic to $S^7$ and they all have the same
isotropy types. \\
	
\end{enumerate}
\end{ex} 

\begin{ex} \label{example514}
Let $q = 5 \cdot 7 = 35$, $q_0 = \frac{35-1}{2} = 17$, $k = 2$, $n = 15$. Here
$A = \{ 1, 2, 3, 4, 6, 8 ,9, 11, 12, 13, 16, 17, 18, 19, 22, 23, 24, 26, 27, 29,
31, 32, 33, 34 \}$, 
$B = \{ 5, 10, 15, 20, 25, 30\}$ and $C = \{ 7, 14, 21, 28 \}$. Here we get isospectral 
non-isometric pairs in three cases:

	\begin{enumerate}
	\item[\underline{Case 1:}]  $q_1 \in A$, $q_2 \in B$, $q_1 \pm q_2 \in A$. 
		So we get, 
		\begin{align*} 
			\psi_{35,2}([q_1, q_2])(z) &= 24z^4  + 6z^3 + 52z^2 +  6z  + 24 \\
			\alpha_{35,2} ([q_1 , q_2])(z) &= 6z^4 + 4z^3 + 8z^2 + 4z + 6 \\ 
			\beta_{35,2} ([q_1 , q_2])(z) &= 4z^4  - 6z^3 + 4z^2 - 6z + 4		
		\end{align*}
	corresponding to 
		\begin{align*}
			L_1 &= L(35: 1, 3, 4, 5, 6, 7, 8, 9, 11, 12, 13, 14, 15, 16, 17) = S^{29} /
G_{1} \\ 
			L_2 &= L(35: 1, 2, 4, 6, 7, 8, 9, 10, 11, 12, 13, 14, 15, 16, 17) = S^{29} /
G_{2}
		\end{align*}
	$L_1$ and $L_2$ are isospectral non-isometric orbifold lens spaces. $G_1 =
\langle g_1 \rangle$ and $G_2 = \langle g_2 \rangle$. 
	\begin{align*}
	\textstyle \sum_{1} &= \big \{ (0, \ldots, 0, x_7, x_8, 0, \ldots, 0, x_{25},
x_{26},0, \ldots, 0) \in S^{29} \text{ with isotropy group} = \langle g_1^7
\rangle \big \} \\ 
	& \bigcup \big \{ (0, \ldots,0,  x_{11}, x_{12}, 0, \ldots, 0, x_{23}, x_{24},
0, \ldots, 0) \in S^{29} \text{ with isotropy group} = \langle g_1^5 \rangle
\big \} \\ \\ 
	\textstyle \sum_{2} &= \big \{ (0, \ldots, 0, x_{15}, x_{16}, 0, \ldots, 0,
x_{25}, x_{26},0, \ldots, 0) \in S^{29} \text{ with isotropy group} = \langle
g_2^7 \rangle \big \} \\ 
	& \bigcup \big \{ (0, \ldots,0,  x_{9}, x_{10}, 0, \ldots, 0, x_{23}, x_{24},
0, \ldots, 0) \in S^{29} \text{ with isotropy group} = \langle g_2^5 \rangle
\big \} 	
	\end{align*}	
	$\sum_1$ and $\sum_2$ are both homeomorphic to $S^{3} \times S^{3}$. 
	
	\item[\underline{Case 2:}] $q_1, q_2 \in A$. 
	
	$(a)$ $q_1 \pm q_2 \in A$. So we get, 
		\begin{align*} 
			\psi_{35,2}([q_1, q_2])(z) &= 24z^4  - 4z^3 + 52z^2 - 4z  + 24 \\
			\alpha_{35,2} ([q_1 , q_2])(z) &= 6z^4 + 4z^3 + 8z^2 + 4z + 6 \\ 
			\beta_{35,2} ([q_1 , q_2])(z) &= 4z^4  + 4z^3 + 4z^2 + 4z + 4		
		\end{align*}	
	corresponding to
		\begin{align*}
			L_3 &= L(35 : 1,3, 5, 6, 7, 8, 9, 10, 11, 12, 13, 14, 15, 16, 17) = S^{29} /
G_{3} \\ 
		\text{ and} \quad L_4 &= L(35: 1, 3, 4, 5, 7, 8, 9, 10, 11, 12, 13, 14, 15,
16, 17) = S^{29} / G_{4}
		\end{align*}
	where $G_3 = \langle g_3 \rangle$ and $G_4 = \langle g_4 \rangle$. Thus, 
		\begin{align*}
		\textstyle \sum_{3} &= \big \{ (0, \ldots, 0, x_5, x_6, 0, \ldots, 0, x_{15},
x_{16}, 0, \ldots, 0, x_{25}, x_{26}, 0, \ldots, 0) \in S^{29} \\ 
		&\qquad \qquad \qquad \qquad \qquad \qquad \qquad \qquad \qquad \text{ with
isotropy group } = \langle g_3^7 \rangle \big \} \\ 
		& \bigcup \big \{ (0, \ldots, 0, x_9, x_{10}, 0 , \ldots, 0, x_{23}, x_{24},
0, \ldots, 0) \in S^{29} \text{ with isotropy group} = \langle g_3^5 \rangle
\big \} \\ \\ 
		\textstyle \sum_{4} &= \big \{ (0, \ldots, 0, x_7, x_8, 0, \ldots, 0, x_{15},
x_{16}, 0, \ldots, 0, x_{25}, x_{26}, 0, \ldots, 0) \in S^{29} \\ 
		&\qquad \qquad \qquad \qquad \qquad \qquad \qquad \qquad \qquad \text{ with
isotropy group } = \langle g_4^7 \rangle \big \} \\ 
		& \bigcup \big \{ (0, \ldots, 0, x_9, x_{10}, 0 , \ldots, 0, x_{23}, x_{24},
0, \ldots, 0) \in S^{29} \text{ with isotropy group} = \langle g_4^5 \rangle
\big \}
		\end{align*}
	$\sum_3$ and $\sum_4$ are homeomorphic to $S^{5} \times S^{3}$. 	
	
	$(b)$ $q_1 + q_2 \in A$, $q_1 - q_2 \in B$. So we get, 
		\begin{align*} 
			\psi_{35,2}([q_1, q_2])(z) &= 24z^4  - 4z^3 + 42z^2 - 4z  + 24 \\
			\alpha_{35,2} ([q_1 , q_2])(z) &= 6z^4 + 4z^3 + 8z^2 + 4z + 6 \\ 
			\beta_{35,2} ([q_1 , q_2])(z) &= 4z^4  + 4z^3 + 14z^2 + 4z + 4		
		\end{align*}		
	corresponding to 
		\begin{align*}
			L_5 &= L(35: 1, 3, 4, 5, 6, 7, 9, 10, 11, 12, 13, 14, 15, 16, 17) = S^{29} /
G_5 \\ 
			\text{ and} \quad L_6 &= L(35: 1, 4, 5, 6, 7, 8, 9, 10, 11, 12, 13, 14, 15,
16, 17) = S^{29} / G_6
		\end{align*}
	where $G_5 = \langle g_5 \rangle$ and $G_6  = \langle g_6 \rangle$. 
		\begin{align*}
			\textstyle \sum_{5} &= \big \{ (0, \ldots, 0, x_7, x_8, 0, \ldots, 0, x_{15},
x_{16}, 0, \ldots, 0, x_{25}, x_{26}, 0, \ldots, 0) \in S^{29} \\ 
			&\qquad \qquad \qquad \qquad \qquad \qquad \qquad \qquad \qquad \text{ with
isotropy group } = \langle g_5^7 \rangle \big \} \\ 
			& \bigcup \big \{ (0, \ldots, 0, x_{11}, x_{12}, 0 , \ldots, 0, x_{23},
x_{24}, 0, \ldots, 0) \in S^{29} \text{ with isotropy group} = \langle g_5^5
\rangle \big \}
		\end{align*}
		
		\begin{align*}
			\textstyle \sum_{6} &= \big \{ (0, \ldots, 0, x_5, x_6, 0, \ldots, 0, x_{15},
x_{16}, 0, \ldots, 0, x_{25}, x_{26}, 0, \ldots, 0) \in S^{29} \\ 
			&\qquad \qquad \qquad \qquad \qquad \qquad \qquad \qquad \qquad \text{ with
isotropy group } = \langle g_6^7 \rangle \big \} \\ 
			& \bigcup \big \{ (0, \ldots, 0, x_9, x_{10}, 0 , \ldots, 0, x_{23}, x_{24},
0, \ldots, 0) \in S^{29} \text{ with isotropy group} = \langle g_6^5 \rangle
\big \}
		\end{align*}
		$\sum_{5}$ and $\sum_{6}$ are homeomorphic to $S^{5} \times S^{3}$. 
		
	\end{enumerate}
\end{ex}

Our final example for the case when $k = 2$ comes when $q$ is even. 

\begin{ex}\label{example515}
Let $q = 2 \cdot 7 = 14$, $q_0 = \frac{14}{2} = 7$, $k =2$ and $n=5$. Here $A =
\{ 1, 3, 5, 9, 11, 13 \}$, $B = \{2, 4, 6, 8, 10, 12 \}$ and $C = \{ 7 \}$. 
Here we get isospctral non-isometric pairs in only one case:

	\begin{enumerate}
	\item[\underline{Case 1:}]  $q_1 \in A$, $q_2 \in B$, $q_1 \pm q_2 \in A$. 
		\begin{align*}
			a_1 &= 2 \sum_{l \in A} \gamma^{l} + 2\sum_{l \in A} \lambda^{l} = 2(1) +
2(-1) = 0 \\ 
			b_1 &= 2 \sum_{l \in B} \gamma^{l} + 2\sum_{l \in B} \lambda^{l} = 2(-1) +
2(-1) = -4 \\
			c_1 &= 2 \sum_{l \in C} \gamma^{l} + 2\sum_{l \in C} \lambda^{l} = 2(-1) +
2(1) = 0
		\end{align*}
		\begin{align*}
			a_2 &= 2 \abs{A} + 4 \sum_{l \in A} \gamma^{l} = 2(6) + 4(1) = 12 + 4 = 16 \\ 
			b_2 &= 2 \abs{B} + 4 \sum_{l \in B} \gamma^{l} = 2(6) + 4(-1) = 12 - 4 = 8 \\ 
			c_2 &= 2 \abs{C} + 4 \sum_{l \in C} \gamma^{l} = 2(1) + 4(-1) = 2 - 4 = -2
		\end{align*}
	So we get, 
		\begin{align*} 
			\psi_{14,2}([q_1, q_2])(z) &= 6z^4 + 16z^2 + 6 \\
			\alpha_{14,2} ([q_1 , q_2])(z) &= 6z^4 + 4z^3 + 8z^2 + 4z + 6 \\ 
			\beta_{14,2} ([q_1 , q_2])(z) &= z^4  - 2z^2 + 1		
		\end{align*}
	corresponding to
		\begin{align*}
			L_1 &= L(14: 1, 2, 4, 5, 7) = S^{9} / G_{1} \\ 
			\text{and} \quad L_2 &= L(14: 1, 4, 5, 6,7) = S^{9} / G_2 
		\end{align*}
	where $G_1 = \langle g_1 \rangle$ and $G_2 = \langle g_2 \rangle$. 
		\begin{align*}
		\textstyle \sum_{1} &= \big \{ (0,0,x_3, x_4, x_5, x_6, 0,0,0,0) \in S^{9}
\text{ with isotropy group} = \langle g_1^7 \rangle \big \} \\
		& \bigcup \big \{ (0, \ldots, 0, x_9, x_{10}) \in S^{9} \text{ with isotropy
group} =  \langle g_1^2 \rangle \big \} \\  \\
		\textstyle \sum_{2} &= \big \{ (0,0,x_3, x_4, 0,0, x_7, x_8, 0,0) \in S^{9}
\text{ with isotropy group} = \langle g_2^7 \rangle \big \} \\
		& \bigcup \big \{ (0, \ldots, 0, x_9, x_{10}) \in S^{9} \text{ with isotropy
group} =  \langle g_2^2 \rangle \big \} 
		\end{align*}
	$\sum_1$ and $\sum_2$ are homeomorphic to $S^{3} \times S^{1}$. \\
	
	\end{enumerate}
\end{ex}

\subsection{Example for $k = 3$}\label{section52}

\begin{ex}\label{example521}
Let $q = 5^2 = 25$, $q_0 = \frac{25 - 1}{2} = 12$, $k = 3$, $n = 9$. Let
$w([p_1, \ldots, p_9 ]) = [q_1, q_2 ,q_3]$. 
Here $A = \{1,2,3,4,6,7,8,9,11, 12, 13, 14, 16, 17, 18, 19, 21, 22, 23, 24 \}$
and \\ $B_1 = \{ 5, 10, 15, 20 \}$. 

We will consider all the possible cases for the various possibilities of
$q_i$'s, $(q_i \pm q_j)$'s and $(q_1 \pm q_2 \pm q_3)$'s lying in $A$ or $B_1$. 
Many of these possibilities will not occur in our present example of $q = 25$.
However, these possibilities are stated because they may occur for higher values
of $q$. 

\begin{enumerate}
\item[\underline{Case 1:}] All the $q_i$'s $\in B_1$. This case does not happen
for $q = 25$ since at most $2$ of the $q_i$'s can be in $B_1$ at one time by the 
definition of $I_{0}(25,3)$. 

\item[\underline{Case 2:}] All the $q_i$'s $\in A$. Since $k = 3$, we can't have
more than $3$ of the $q_i \pm q_j$ $(1 \leq i < j \leq 3 )$ belonging to $B_1$.
Also, at 
most only one of the $q_1 \pm q_2 \pm q_3$ can be $\equiv 0 \pmod{25}$. 

Further, if one of the $q_1 \pm q_2 \pm q_3$ is congruent to $0 \pmod{25}$, then
we can't have any other of the  $q_1 \pm q_2 \pm q_3$ belong to $B_1$.  
Also, at most $1$ of the $q_1 \pm q_2 \pm q_3$ can be in $B_1$. 

Further, if one of the $q_1 \pm q_2 \pm q_3$ is congruent to $0 \pmod{25}$, 
then at most $1$ of the $q_i \pm q_j$ $(\text{for } 1 \leq i < j
\leq 3)$ can be in $B_1$. Similarly, if one of the $q_1 \pm q_2
\pm q_3 \in B_1$, then at most $1$ of the $q_i \pm q_j$ can be in $B_1$. 
We note that the above results hold true for all $q = P^2$, where $P$ is any odd
prime. 

We now look at the various sub-cases for Case 2. 

$(a)$ All of the $q_1 \pm q_2 \pm q_3 \in A$ and all of the $q_i \pm q_j \in A$
$(1 \leq i < j \leq 3)$. This case does not occur for $q = 25$. 

$(b)$ All of the $q_1 \pm q_2 \pm q_3 \in A$ and exactly one of the $q_i \pm q_j
\in B_1$ $(1 \leq i < j \leq 3)$. This case does not occur for $q = 25$. 

$(c)$ All of the $q_1 \pm q_2 \pm q_3 \in A$ and exactly two of the $q_i \pm q_j
\in B_1$ $(1 \leq i < j \leq 3)$. Again, this case does not occur for $q = 25$. 

$(d)$ exactly one of the  $q_1 \pm q_2 \pm q_3 \in B_1$ and all of the $q_i \pm
q_j \in A$ $(1 \leq i < j \leq 3)$. This case does not occur for $q = 25$. 

$(e)$ exactly one of the $q_1 \pm q_2 \pm q_3 \equiv 0 \pmod{q}$ and all of $q_i
\pm q_j \in A$ $(1 \leq i < j \leq 3)$. This case also does not occur for $q=
25$. 

$(f)$ All of the $q_1 \pm q_2 \pm q_3 \in A$ and exactly $3$ of the $q_i \pm q_j
\in B_1$ $(1 \leq i < j \leq 3)$.  In this case, we get isospectral, non-isometric 
pairs since we get,

\begin{align*}
	\psi_{25,3} ([ q_1, q_2 ,q_3])(z) &= 20z^6 + 30z^4 + 30z^2 + 20 \\
	\alpha_{25,3}^{(1)} ([q_1, q_2, q_3 ])(z) &= 4z^6 + 6z^5 + 30z^4 + 20z^3 + 30z^2 +
6z + 4 
\end{align*}

corresponding to
\begin{align*}
	L_1 &= L(25: 1, 4, 5, 6, 7, 8, 9, 10 ,11) = S^{17} / G_{1} \\ 
	\text{and} \quad L_2 &= L(25: 1, 2, 3, 5, 7, 8, 9, 10 , 12) = S^{17} / G_2 
\end{align*}

where $G_1 = \langle g_1 \rangle$ and $G_2 = \langle g_2 \rangle$. 
\begin{align*}
\textstyle \sum_{1} &= \big \{ (0, \ldots, 0, x_5, x_6, 0, \ldots, 0, x_{15},
x_{16}, 0, 0) \in S^{17} \text{ with isotropy group}= \langle g_1^5 \rangle \big
\} \\ 
\textstyle \sum_{2} &= \big \{ (0, \ldots, 0, x_7, x_8 , 0, \ldots, 0, x_{15},
x_{16}, 0 , 0) \in S^{17} \text{ with isotropy group} = \langle g_2^5 \rangle
\big \}
\end{align*}

$(h)$ exactly one of the $q_1 \pm q_2 \pm q_3$ is congruent to $0 \pmod{25}$ and
one of the $q_i \pm q_j$ $(\text{for } 1 \leq i < j \leq 3)$ is in $B_1$. We
again get isospectral non-isometric pairs here since,

\begin{align*}
	\psi_{25,3} ([ q_1, q_2 ,q_3])(z) &= 20z^6 + 50z^4 +10z^3+  50z^2 + 20 \\
	\alpha_{25,3} ([q_1, q_2, q_3 ])(z) &= 4z^6 + 6z^5 + 10z^4 + 10z^3 + 10z^2 + 6z +
4 	
\end{align*}

corresponding to 
\begin{alignat*}2
	L_3 &= L(25: 1, 2,3, 4, 5, 6, 7, 8, 10) = S^{17} / G_3, &\quad G_3 &= \langle
g_3 \rangle \\ 
	L_4 &= L(25: 1, 2, 3, 4, 5, 6, 7, 9, 10) = S^{17} / G_4, &\quad G_4 &= \langle
g_4 \rangle \\ 
	L_5 &= L(25: 1, 2,3, 4, 5, 6, 8,9, 10) = S^{17} / G_5, &\quad G_5 &= \langle
g_5 \rangle \\ 
	L_6 &= L(25: 1,2,3,4,5,7,8,9, 10) = S^{17} / G_6, &\quad G_6 &= \langle g_6
\rangle \\ 
	L_7 &= L(25: 1,2,4,5,6,7,8,9,10) = S^{17} / G_7, &\quad G_7 &= \langle g_7
\rangle \\ 
	L_8 &= L(25: 1, 4,5,6,7,8,9,10,12) = S^{17} / G_8, &\quad G_8 &= \langle g_8
\rangle \\ 
	L_9 &= L(25: 1, 3, 5, 6, 7, 9, 10, 11, 12) = S^{17} / G_9, &\quad G_9 &=
\langle g_9 \rangle \\ 
	L_{10} &= L(25: 1, 2, 5, 6, 7, 8, 9, 10, 12) = S^{17} / G_{10}, &\quad G_{10}
&= \langle g_{10} \rangle
\end{alignat*}

So, in this case we get a family of $8$ orbifold lens spaces that are
isospectral but mutually non-isometric. 
\begin{align*}
\textstyle \sum_{3} &= \big \{ (0, \ldots, 0, x_9, x_{10}, 0, \ldots, 0, x_{17},
x_{18} ) \in S^{17} \text{ with isotropy group}=\langle g_3^5 \rangle \big \} \\
\textstyle \sum_{4} &= \big \{ (0, \ldots, 0, x_9, x_{10}, 0, \ldots, 0, x_{17},
x_{18} ) \in S^{17} \text{ with isotropy group}=\langle g_4^5 \rangle \big \} \\
\textstyle \sum_{5} &= \big \{ (0, \ldots, 0, x_9, x_{10}, 0, \ldots, 0, x_{17},
x_{18} ) \in S^{17} \text{ with isotropy group}=\langle g_5^5 \rangle \big \} \\
\textstyle \sum_{6} &= \big \{ (0, \ldots, 0, x_9, x_{10}, 0, \ldots, 0, x_{17},
x_{18} ) \in S^{17} \text{ with isotropy group}=\langle g_6^5 \rangle \big \} 
\end{align*}
\begin{align*}
\textstyle \sum_{7} &= \big \{ (0, \ldots, 0, x_7, x_{8}, 0, \ldots, 0, x_{17},
x_{18} ) \in S^{17} \text{ with isotropy group}=\langle g_7^5 \rangle \big \} \\
\textstyle \sum_{8} &= \big \{ (0, \ldots, 0, x_5, x_{6}, 0, \ldots, 0, x_{15},
x_{16},0,0 ) \in S^{17} \text{ with isotropy group}=\langle g_8^5 \rangle \big
\} \\
\textstyle \sum_{9} &= \big \{ (0, \ldots, 0, x_5, x_{6}, 0, \ldots, 0, x_{13},
x_{14},0, \ldots,0 ) \in S^{17} \text{ with isotropy group}=\langle g_9^5
\rangle \big \}  \\
\textstyle \sum_{10} &= \big \{ (0, \ldots, 0, x_5, x_{6}, 0, \ldots, 0, x_{15},
x_{16},0,0 ) \in S^{17} \text{ with isotropy group}=\langle g_{10}^5 \rangle
\big \}
\end{align*}
All of the $\sum_i$ $(\text{for } i = 3,4, \ldots, 10)$ are homeomorphic to
$S^{3}$. 

$(h)$ exactly one of the $q_1 \pm q_2 \pm q_3$ is congruent to $0 \pmod{25}$ and
one of the $q_i \pm q_j$ $(\text{for } 1 \leq i < j \leq 3)$ is in $B_1$. 
Here we get, 
\begin{align*}
	\psi_{25,3} ([ q_1, q_2 ,q_3])(z) &= 20z^6 + 50z^4 - 40z^3+  50z^2 + 20 \\
	\alpha_{25,3} ([q_1, q_2, q_3 ])(z) &= 4z^6 + 6z^5 + 10z^4 + 10z^3 + 10z^2 + 6z +
4 	
\end{align*}

corresponding to 
\begin{alignat*}2
	L_{11} &= L(25: 1,3,4,5,6,7,8,9,10) = S^{17} / G_{11}, & \quad G_{11} &=
\langle g_{11} \rangle \\
	L_{12} &= L(25: 1,3,5,7, 8, 9, 10, 11, 12) = S^{17} / G_{11}, & \quad G_{12} &=
\langle g_{12} \rangle 
\end{alignat*}
\begin{align*}
\textstyle \sum_{11} &= \big \{ (0, \ldots, 0, x_7, x_{8}, 0, \ldots, 0, x_{17},
x_{18} ) \in S^{17} \text{ with isotropy group}=\langle g_{11}^5 \rangle \big \}
\\
\textstyle \sum_{12} &= \big \{ (0, \ldots, 0, x_5, x_{6}, 0, \ldots, 0, x_{13},
x_{13},0, \ldots,0 ) \in S^{17} \text{ with isotropy group}=\langle g_{12}^5
\rangle \big \}
\end{align*}
$\sum_{11}$ and $\sum_{12}$ are homeomorphic to $S^{3}$. \\

$2(a)$ - $2(h)$ are all of the possible cases when all the $q_i$'s $\in A$. 

\item[\underline{Case 3:}] Two of the $q_i$'s $\in A$ and one of the $q_i$'s
$\in B_1$. In this case we will have at most one of the $q_1 \pm q_2 \pm q_3$ 
congruent to $0 \pmod{25}$. 

Also, we can have at most one of the $q_i \pm q_j$ $(1 \leq i < j \leq 3)$ in
$B_1$. 

Further, it can be shown that at most two of the $q_1 \pm q_2 \pm
q_3$ can be in $B_1$. Also, if one of the $q_1 \pm q_2 \pm q_3$ is congruent
to $0 \pmod{25}$, then at most one of the remaining $q_1 \pm q_2 \pm q_3$ can
belong to $B_1$. In fact, it can be shown that exactly one of the 
remaining $q_1 \pm q_2 \pm q_3$ must belong to $B_1$. 

All of these results can be shown to be true $q = P^2$, where $P$ is any odd
prime. Now we consider all the sub-cases for Case 3. 

$(a)$ If all the $q_1 \pm q_2 \pm q_3$ belong to $A$ and exactly one of the $q_i
\pm q_j$ $(1 \leq i < j \leq 3)$ belongs to $B_1$. This case does not 
occur for $q = 25$. 

$(b)$ exactly one of the $q_1 \pm q_2 \pm q_3$ belongs to $B_1$ and the
remaining belong to $A$. This case does not occur for $q = 25$. 

$(c)$ If all of the $q_1 \pm q_2 \pm q_3$ belong to $A$ and all of the $q_i \pm
q_j$ $(1 \leq i < j \leq 3)$ belong to $A$. Then we get
	 
\begin{align*}
	\psi_{25,3} ([ q_1, q_2 ,q_3])(z) &= 20z^6 + 10z^5 + 60z^4 + 20z^3 + 60z^2 + 10z +
20 \\
	\alpha_{25,3}^{(1)} ([q_1, q_2, q_3 ])(z) &= 4z^6 - 4z^5 - 4z  + 4
\end{align*}

corresponding to 
\begin{alignat*}2
	L_{13} &= L(25: 1,2,3,4,5,6,7,8,9) = S^{17} / G_{13}, &\quad G_{13} &= \langle
g_{13} \rangle \\ 
	L_{14} &= L(25: 1,2,3,4,5,6,7,8,11) = S^{17} / G_{14}, &\quad G_{14} &= \langle
g_{14} \rangle \\ 
	L_{15} &= L(25: 1,2,3,4,5,6,7,9,12) = S^{17} / G_{15}, &\quad G_{15} &= \langle
g_{15} \rangle
\end{alignat*}

\begin{alignat*}2
	L_{16} &= L(25: 1,3,5,6,7,8,9,11,12) = S^{17} / G_{16}, &\quad G_{16} &=
\langle g_{16} \rangle \\ 
	L_{17} &= L(25: 1,2,4,5,6,7,8,9,12) = S^{17} / G_{17}, &\quad G_{17} &= \langle
g_{17} \rangle 
\end{alignat*}
	
We get a family of $5$ orbifold lens spaces that are non-isometric and
isospectral. 
\begin{align*}
\textstyle \sum_{13} &= \big \{ (0, \ldots, 0, x_9, x_{10}, 0, \ldots, 0 ) \in
S^{17} \text{ with isotropy group}=\langle g_{13}^5 \rangle \big \} \\
\textstyle \sum_{14} &= \big \{ (0, \ldots, 0, x_9, x_{10}, 0, \ldots, 0) \in
S^{17} \text{ with isotropy group}=\langle g_{14}^5 \rangle \big \} \\
\textstyle \sum_{15} &= \big \{ (0, \ldots, 0, x_9, x_{10}, 0, \ldots, 0) \in
S^{17} \text{ with isotropy group}=\langle g_{15}^5 \rangle \big \} \\
\textstyle \sum_{16} &= \big \{ (0, \ldots, 0, x_5, x_{6}, 0, \ldots, 0) \in
S^{17} \text{ with isotropy group}=\langle g_{16}^5 \rangle \big \} \\
\textstyle \sum_{17} &= \big \{ (0, \ldots, 0, x_7, x_{8}, 0, \ldots, 0) \in
S^{17} \text{ with isotropy group}=\langle g_{17}^5 \rangle \big \} 
\end{align*}
All the $\sum_{i}$'s $(i = 13, \ldots, 17)$ are homeomorphic to $S^{1}$. 

$(d)$ exactly two of the $q_1 \pm q_2 \pm q_3$ belong to $B_1$ and exactly one
of the $q_i \pm q_j$ $(1 \leq i < j \leq 3)$ belongs to $B_1$. 

Here we get, 
\begin{align*}
	\psi_{25,3} ([ q_1, q_2 ,q_3])(z) &= 20z^6 + 10z^5 + 50z^4 + 40z^3 + 50z^2 + 10z +
20 \\
	\alpha_{25,3}^{(1)} ([q_1, q_2, q_3 ])(z) &= 4z^6 - 4z^5 + 10z^4 - 20z^3 + 10z^2 -
4z  + 4
\end{align*}

corresponding to 
\begin{alignat*}2
	L_{18} &= L(25: 1,2,3,4,5,6,7,8,12) = S^{17} / G_{18}, &\quad G_{18} &= \langle
g_{18} \rangle \\ 
	L_{19} &= L(25: 1, 2, 3, 4, 5, 6, 7, 9, 11) = S^{17} / G_{19}, &\quad G_{19} &=
\langle g_{19} \rangle 
\end{alignat*}

We have
\begin{align*}
\textstyle \sum_{18} &= \big \{ (0, \ldots, 0, x_9, x_{10}, 0, \ldots, 0 ) \in
S^{17} \text{ with isotropy group}=\langle g_{18}^5 \rangle \big \} \\
\textstyle \sum_{19} &= \big \{ (0, \ldots, 0, x_9, x_{10}, 0, \ldots, 0) \in
S^{17} \text{ with isotropy group}=\langle g_{19}^5 \rangle \big \} \\
\end{align*}

$(e)$ One of the $q_1 \pm q_2 \pm q_3$ is congruent to $0 \pmod{25}$, and one
one of the $q_1 \pm q_2 \pm q_3$ is in $B_1$, and exactly one of the
$q_i \pm q_j$ $(1 \leq i < j \leq 3)$ is in $B_1$. 

Here we get, 
\begin{align*}
	\psi_{25,3} ([ q_1, q_2 ,q_3])(z) &= 20z^6 + 10z^5 + 50z^4 - 10z^3 + 50z^2 + 10z +
20 \\
	\alpha_{25,3}^{(1)} ([q_1, q_2, q_3 ])(z) &= 4z^6 - 4z^5 + 10z^4 - 20z^3 + 10z^2 -
4z  + 4
\end{align*}

corresponding to 
\begin{alignat*}2
	L_{20} &= L(25: 1,2, 3, 4, 6, 7, 8, 10, 12) = S^{17} / G_{20}, &\quad G_{20} &=
\langle g_{20} \rangle \\ 
	L_{21} &= L(25: 1, 4, 6, 7, 8, 9, 10, 11 ,12) = S^{17} / G_{21}, &\quad G_{21}
&= \langle g_{21} \rangle \\ 
	L_{22} &= L(25: 1, 3, 4, 5, 6, 7, 9 , 11, 12) = S^{17} / G_{22}, &\quad G_{22}
&= \langle g_{22} \rangle 
\end{alignat*}
We have
\begin{align*}
\textstyle \sum_{20} &= \big \{ (0, \ldots, 0, x_{15}, x_{16}, 0, 0 ) \in S^{17}
\text{ with isotropy group}=\langle g_{20}^5 \rangle \big \} \\
\end{align*}
\begin{align*}
\textstyle \sum_{21} &= \big \{ (0, \ldots, 0, x_{13}, x_{14}, 0, \ldots, 0) \in
S^{17} \text{ with isotropy group}=\langle g_{21}^5 \rangle \big \} \\
\textstyle \sum_{22} &= \big \{ (0, \ldots, 0, x_7, x_{8}, 0, \ldots, 0) \in
S^{17} \text{ with isotropy group}=\langle g_{22}^5 \rangle \big \} 
\end{align*}

There are no other sub-cases for Case 3. 

\item[\underline{Case 4:}]  One of the $q_i$'s $\in A$ and two of the $q_i$'s
$\in B_1$. In this case, $q_1 \pm q_2 \pm q_3$ will always belong to $A$, 
and exactly two of the $q_i \pm q_i$ $(1 \leq i < j \leq 3)$ will belong to
$B_1$. There are no other variations that will occur in this case. 

Here we get, 

\begin{align*}
	\psi_{25,3} ([ q_1, q_2 ,q_3])(z) &= 20z^6 + 20z^5 + 40z^4 + 40z^3 + 40z^2 + 20z +
20 \\
	\alpha_{25,3}^{(1)} ([q_1, q_2, q_3 ])(z) &= 4z^6 - 14z^5 + 20z^4 - 20z^3 + 20z^2
- 14z  + 4
\end{align*}

corresponding to $L(25: 1, 2, 3, 4, 6, 7, 8, 9, 11)$ and we do not get
isospectral pairs. Note that this lens space is a manifold. 
\end{enumerate}
\end{ex}

As this example illustrates, we can extend our technique for $k=2$ to higher
values of $k$ and we will get many examples of isospectral non-isometric 
orbifold lens spaces. At the same time, the example also illustrates the
difficulty in accounting for all the possible cases as the value of $k$ is
increased.






\addtocontents{toc}{\protect\vspace{2ex}}


\begin{thebibliography}{xx}

\bibitem[BCDS]{BCDS} P. Buser, J. Conway, P. Doyle and K. Semmler, {\it Some planar isospectral domains},
Internat. Math. Res. Notices. {\bf 9} (1994), 391ff., approx. 9 pp. (electronic).

\bibitem[BGM]{BGM} M. Berger, P. Gaudachon and E. Mazet, {\it Le spectre d'une
vari\'et\'e riemannienne}, Lecture notes in Mathematics 194, 
Springer-Verlag, Berlin-Heidelberg-New York, 1971. 

\bibitem[BW]{BW} P. B\'erard and D. Webb, {\it On ne peut pas entendre l\'orientabilit\'e d\'une surface}, 
C. R. Acad. Sci. Paris S´er. I Math. {\bf 320} (1995), no. 5, 533–536.

\bibitem[Chi]{Chi} Chiang, Yuan-Jen, {\it Spectral Geometry of V-Manifolds and
its Application to Harmonic Maps}, Proc. Symp. Pure Math. {\bf 54} part 1
(1993), 93--99. 

\bibitem[DGGW]{DGGW} E. Dryden, C. Gordon, S. Greenwald and D. Webb, 
{\it Asymptotic expansion of the heat kernel for orbifolds}, Michigan Math J. 56 (2008), 205--238.

\bibitem[DR]{DR} P. Doyle and J. Rossetti, {\it Isospectral hyperbolic surfaces having matching geodesics},
preprint, ArXiv math.DG/0605765.

\bibitem[GR]{GR} C. S. Gordon and J. Rossetti, {\it Boundary volume and length spectra of Riemannian
manifolds: what the middle degree Hodge spectrum doesn't reveal}, Ann. Inst. Fourier, {\bf 53}
(2003), no. 7, 2297–2314.

\bibitem[Gi]{Gi} P. B. Gilkey, {\it On spherical space forms with meta-cyclic fundamental group which 
are isospectral but not equivariant cobordant}. Compositio Mathematica, {\bf 56} no. 2 (1985), p. 171-200 

\bibitem[GoM]{GoM} R. Gornet and J. McGowan, {\it Lens spaces, isospectral on forms but not on functions}, 
London Math. Soc. J. of Computation {\bf 9} (2006)  270-286.
    
\bibitem[I1]{I1} A. Ikeda, {\it On lens spaces which are isospectral but not
isometric}, Ann. scient. \'Ec. Norm. Sup. $4^{e}$ s\'eries, t. 13, 303--315.

\bibitem[I2]{I2} A. Ikeda, {\it On the spectrum of a riemannian manifold of
positive constant curvature}, Osaka J. Math., {\it 17} (1980), 75--93.

\bibitem[IY]{IY} A. Ikeda and Y. Yamamoto, {\it On the spectra of a
3-dimensional lens space}, Osaka J. Math., {\it 16} (1979), 447--469.

\bibitem[K]{K} M. Kac, {\it Can one hear the shape of a drum?}, Amer. Math.
Monthly 73 (1966) no. 4m Part II, 1--23.

\bibitem[M]{M} J. Milnor, {\it Eigenvalues of the Laplace operator on certain manifolds}, 
Proc. Nat. Acad. Sci. USA {\bf 51} (1964), 542.

\bibitem[PS]{PS} E. Proctor and E. Stanhope, {\it An Isospectral Deformation on an
Orbifold Quotient of a Nilmanifold}, Preprint, ArXiv math. 0811.0794

\bibitem[RSW]{RSW} J. Rossetti, D. Schueth and M. Weilandt, {\it Isospectral orbifolds 
with different maximal isotropy orders}, Ann. Glob. Anal. Geom. {\bf 34} (2008), 351 - 366

\bibitem[Sat]{Sat} I. Satake, {\it On a generalization of the notion of manifold}, 
Proc. Nat. Acad. Sci. U.S.A. {\bf 42} (1956), 359--363.

\bibitem[Sc]{Sc} P. Scott, {\it The geometries of 3-manifolds} Bull.London Math.Soc.
{\bf 15} (1983), no. 5, 401--487.

\bibitem[SSW]{SSW} N. Shams, E. Stanhope, and D. Webb, {\it One Cannot Hear
Orbifold Isotropy Type}, Archiv der Math (Basel) 87 (2006), no.4, 375-384.

\bibitem[S1]{S1} E. Stanhope, {\it Hearing Orbifold Topology}, Ph.D. Thesis,
Dartmouth College, 2002. 

\bibitem[S2]{S2} E. Stanhope, {\it Spectral bounds on orbifold isotropy}, Annals
of Global Analysis and Geometry {\bf 27} (2005), no. 4, 355--375.

\bibitem[V]{V} M. F. Vign\'eras, {\it Vari\'et\'es Riemanniennes isospectrales et non isom´etriques}, 
Ann. of Math. {\bf 112} (1980), 21–32.

\end{thebibliography}
\end{document}